\documentclass[12pt]{article}
\usepackage{amssymb,amsmath,amsfonts,amsthm,enumerate,color}
\usepackage[toc,page,title,titletoc,header]{appendix} 
\usepackage{graphicx}
\usepackage{indentfirst}
\usepackage{multicol}
\usepackage{booktabs}

\numberwithin{equation}{section}

\numberwithin{equation}{section}

\def \Vh0{\stackrel{\circ}{V}_h} \def\to{\rightarrow}

\newcommand{\q}{\quad}

\def\ms{\medskip}  \def\ss{\smallskip}

  \def\x{{\boldsymbol x}}

\newcommand{\lc}
{\mathrel{\raise2pt\hbox{${\mathop<\limits_{\raise1pt\hbox
{\mbox{$\sim$}}}}$}}}

\newcommand{\gc}
{\mathrel{\raise2pt\hbox{${\mathop>\limits_{\raise1pt\hbox{\mbox{$\sim$}}}}$}}}

\newcommand{\ec}
{\mathrel{\raise2pt\hbox{${\mathop=\limits_{\raise1pt\hbox{\mbox{$\sim$}}}}$}}}

\def\bb{\begin{equation}} \def\ee{\end{equation}}

\def\beqn{\begin{eqnarray}}  \def\eqn{\end{eqnarray}}

\def\beqnx{\begin{eqnarray*}} \def\eqnx{\end{eqnarray*}}

\def\bn{\begin{enumerate}} \def\en{\end{enumerate}}

\def\bd{\begin{description}} \def\ed{\end{description}}

\makeatletter

\newenvironment{figurehere}
  {\def\@captype{figure}}
  {}
\makeatother

\title{A Multilevel Sampling Algorithm \\ 
for Locating Inhomogeneous Media}

\begin{document}

\author{Keji Liu\thanks{Department of Mathematics,
The Chinese University of Hong Kong, Shatin, Hong Kong. ({\tt
kjliu@math.cuhk.edu.hk})} ~~\quad Jun Zou\thanks{Department of
Mathematics, The Chinese University of Hong Kong, Shatin, Hong Kong.
The work of this author  was substantially supported by Hong Kong RGC grants
(Projects 405110 and 404611).
({\tt zou@math.cuhk.edu.hk}) 
}}


\maketitle

\begin{abstract}
In the reconstruction process of
unknown multiple scattering objects in inverse medium scattering problems,
the first important step is to effectively locate some approximate domains that contain
all inhomogeneous media. Without such an effective step, one may have to take a much larger 
computational domain than actually needed in the reconstruction of all 
scattering objects, thus resulting in a huge additional computational efforts. 
In this work we propose a simple and efficient multilevel reconstruction algorithm to help locate 
an accurate position and shape  of each inhomogeneous medium. 
Then other existing effective but computationally 
more demanding reconstruction algorithms may be applied in these initially located computational domains
to achieve more accurate shapes of the scatter and the contrast values over each medium domain.  
The new algorithm exhibits several strengths: robustness against noise, requiring less incidences, 
fast convergence, flexibility to deal with scatterers of special shapes, and advantages 
in computational complexity.
\end{abstract}



\section{Introduction}
This works investigates the numerical identification of inhomogeneous medium scatterers 
by scattered fields. The inverse scattering problem 
can find wide applications in medicine, geophysics, biological studies. 
A large variety of numerical reconstruction methods are available
    in literature, such as
    the time-reversal multiple signal classification (MUSIC) method \cite{Kir,MarGruSim}, the contrast source inversion (CSI) method \cite{BerKle,BerBro,AbuBer}, the continuation method \cite{bao05}, 
    the subspace-based optimization method (SBOM) \cite{Chen09,Chen10}, the linear sampling or 
    probing methods (LSM) \cite{ColKre,Pot,LLZ08},  
    the parallel radial bisection method (PRBM) \cite{LiuXuZou}, etc. 
In order to carry out any of these methods for the reconstruction of
unknown multiple scattering objects, 
the first important step is to effectively locate some approximate domains that contain
all scattering objects. Without such an effective step, one may have to take a much larger 
computational domain than actually needed in the reconstruction of all 
scattering objects.
In particular, when multiple separated objects are present, and at least two of them are
far away from each other, then one may need to set an initial sampling or computational domain 
to be much larger in order to ensure a safe covering of all such objects. 
One may easily have 
an initial computational domain with an area or a volume of 30 or 40 times as large as 
the actual region required to cover all unknown inhomogeneous medium objects. 
A much larger computational domain results usually in a huge additional computational effort
for the entire numerical reconstruction process, considering the severe ill-posedness and 
strong nonlinearity of inverse medium scattering problems. 

So it is of great significance for the reconstruction process of an inverse medium  problem 
to have an effective step that helps locate the initial regions that cover each of the scattering object. 
In addition, this first step should be less expensive computationally and easy to implement 
numerically. 
It is mostly challenging to realize this target, and to provide an acceptable initial location of 
each scattering object at the same time.
A direct sampling method was proposed recently in \cite{Jin12} for the purpose. 
The algorithm is computationally very cheap 
as it involves computing only the inner product of the scattered field with fundamental solutions 
located at sampling points. 
    In this work, we will propose a new algorithm for the purpose, and it is 
    completely different from the one in  \cite{Jin12}. 
    This new algorithm is an iterative one, also very cheap; only three matrix-vector 
    multiplications are needed at each iteration, without any matrix inversion or solutions 
    of linear systems involved. Most interestingly, the algorithm can first separate 
    all disjoint objects quickly, usually in a few iterations, 
    then refine its approximation successively and finally provide a good approximate domain 
    for each separate object. 
    
    It is worth mentioning that the multilevel algorithm to be presented here is 
    essentially different in nature from the multilevel linear
    sampling method developed in \cite{LLZ08}: 
    the new method is much less sensitive to the so-called cut-off values, it 
    works with much less incident fields, and it does not need to solve ill-posed 
    far-field equations at sampling points. 
    In addition, the new algorithm is robust in the presence of noise and less sensitive
    to the noise level. 
    More importantly, unlike most existing methods,
    the new method does not involve any optimization process or matrix inversions, 
    so it can be viewed as a direct sampling method.  
    Another nice feature of the new algorithm is that it is self-adaptive, that is, at each iteration  
    it can remedy the possible errors from the previous iterations. 
    With an effective initial location of each object,
    we may then apply any existing efficient but computationally more demanding 
    methods, e.g., the methods in \cite{bao05}  \cite{Chen10} \cite{BerBro}, 
    for further refinement of the estimated shape of each 
scattering object as well as for recovery of the contrast profiles of different media.

\section{Problem description}

Consider an inverse scattering problem where 
the scatterer $\Omega$, 
possibly consisting of several separated disjoint components, is located in a homogeneous background
medium 
$\mathbb{R}^d$ ($d=2,3$). 
We assume that the scattered obstacles are illuminated 
successively by a number of plane wave incident fields 
$u^{inc}_j(\boldsymbol x)$, $j=1,2,\cdots,N_i$.
For each plane wave incidence, the scattered field $u^{sca}(\boldsymbol x^s_q)$ is measured by the receivers at locations $\boldsymbol x^s_1$, $\cdots$, $\boldsymbol x^s_{N_s}$; see Figure~\ref{fig:example}
for the incidences and receivers located on a circle $\mathrm{S}$.
     \begin{figurehere}
     \begin{center}
     \vskip -0.1truecm
           \scalebox{0.4}{\includegraphics{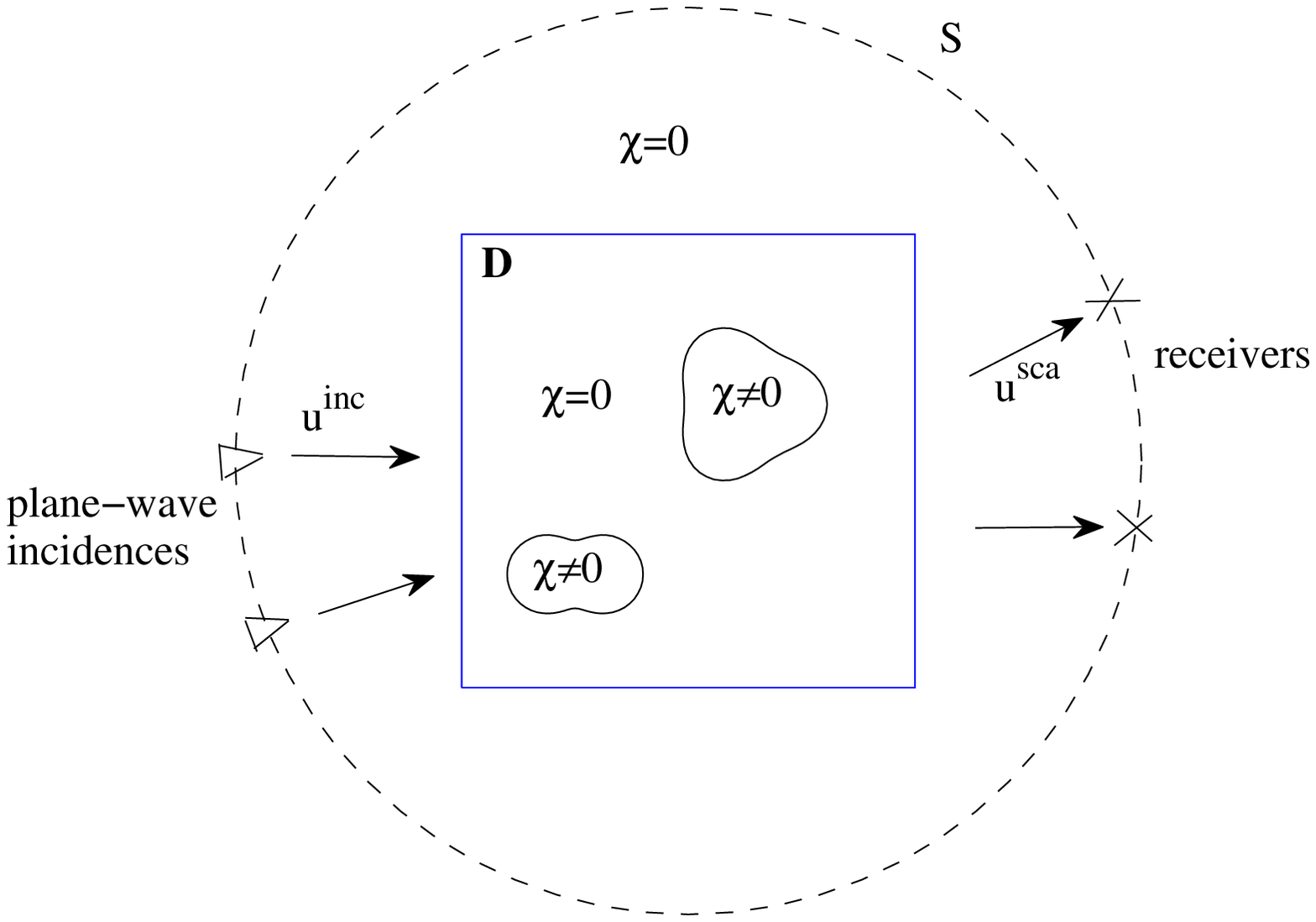}}\\
            \vskip -0.4truecm
     \caption{\small Geometrical model of the scattering problem.}
     \label{fig:example}
     \end{center}
     \end{figurehere}
The inverse scattering problem is to determine the contrast function or index of refraction, 
$\chi(\boldsymbol x)$ for any point $\boldsymbol x$ varying in the scatterer $\Omega$, 
given a set of $N_iN_s$ scattering data $u^{sca}$. 
The contrast $\chi$ has a very important property, namely it vanishes outside the scattering obstacles. 
For each incident field, the total field  $u_j$  satisfies the Helmholtz equation \cite{ColKre}:
    \begin{equation}\label{eq:Etot}
    \Delta u_j(\boldsymbol x)+k^2(\chi(\boldsymbol x)+1)u_j(\boldsymbol x)=0,\quad \boldsymbol x\in \mathrm{\mathbb{R}^d},
    \end{equation}
    where $k$ is the wavenumber of the homogeneous background media. 
    The total field  $u_j$  can be represented by the intergral equation \cite{ColKre},
    \begin{equation}\label{eq:EtotExp1}
    u_j(\boldsymbol x)=u^{inc}_j(\boldsymbol x)+k^2\int_{\Omega}g(\boldsymbol x,\boldsymbol x')\chi(\boldsymbol x') u_j(\boldsymbol x')dv(\boldsymbol x'),
    \end{equation}
    where $g(\boldsymbol x,\boldsymbol x')$ is the Green's function of the homogeneous background
    medium: 
\[
   g(\boldsymbol x,\boldsymbol x')=
    \left\{\begin{array}{cll}
    \frac{i}{4}H^{(1)}_0(k|\boldsymbol x-\boldsymbol x'|) & \textrm{for} & d=2,\\
    \frac{e^{ik|\boldsymbol x-\mathbf{\boldsymbol x'}|}}{4\pi|\boldsymbol x-\mathbf{\boldsymbol x'}|}
    & \textrm{for} & d=3\\
    \end{array}
    \right.
    \]
where $H_0^{(1)}$ is the zero-order Hankel function of first kind. We note that 
the total field $u$ may stand for the acoustic pressure in an acoustic scattering problem, 
or for the electric field vector in an electromagnetic scattering, or for the particle-velocity vector in 
an elastodynamic scattering. 
The scattered field is measured 
    on the boundary of a domain $\mathrm{S}$, which is sitting outside the scatterer $\Omega$. 
    As the scatterer $\Omega$ is unknown, we shall introduce a sampling 
    domain $\mathbf{D}$ that completely cover the scatterer $\Omega$. 
    As the contrast function $\chi$ vanishes outside $\Omega$, 
    with the help of (\ref{eq:EtotExp1})  we can write the scattered field as 
    \begin{eqnarray}\label{eq:Esca1}
    u^{sca}_j(\boldsymbol x)=u_j(\boldsymbol x)-u^{inc}_j(\boldsymbol x)
                                 =k^2\int_\mathrm{D}g(\boldsymbol x,\boldsymbol x')\chi(\boldsymbol x')u_j(\boldsymbol x')dv(\boldsymbol x'),\quad \boldsymbol x\in \mathrm{S}.
    \end{eqnarray}
    
For sake of convenience, we shall often introduce the contrast source function 
    \begin{equation}\label{eq:I}
    w_j(\boldsymbol x)=\chi(\boldsymbol x)u_j(\boldsymbol x),\quad \boldsymbol x\in \mathrm{D}\,.
    \end{equation}
Then we can write (\ref{eq:EtotExp1}) and (\ref{eq:Esca1}) in the following more compact forms
    \begin{equation}\label{eq:IMat}
    w_j(\boldsymbol x)=\chi(\boldsymbol x)u^{inc}_j(\boldsymbol x)+\chi(\boldsymbol x)
    (G_Dw_j)(\boldsymbol x),\quad \boldsymbol x\in\mathrm{D},
    \end{equation}
    and 
    \begin{equation}\label{eq:EsacMat}
    u^{sca}_j(\boldsymbol x)=(G_Sw_j)(\boldsymbol x),\quad \boldsymbol x\in\mathrm{S},  
    \end{equation}
    where $G_D$ and $G_S$ are two integral operators given by 
    \beqnx 
    (G_Dw)(\boldsymbol x)&=&k^2\int_Dg(\boldsymbol x,\boldsymbol x')w(\boldsymbol x')dv(\boldsymbol x')
    \q \forall\,\boldsymbol x\in\mathrm{D}\,, \\
    (G_Sw)(\boldsymbol x)&=&k^2\int_Dg(\boldsymbol x,\boldsymbol x')w(\boldsymbol x')dv(\boldsymbol x')
    \q \forall\,\boldsymbol x\in\mathrm{S}. 
   \eqnx 
    Equations (\ref{eq:IMat}) and (\ref{eq:EsacMat})  will be two fundamental equations 
    for our proposed multilevel initialization algorithm.

\section{Approximate contrast source by backpropagation}

We can easily see that the support of the contrast source function 
$w=\chi u$ describes the exact locations and geometries of all the inhomogeneous 
media, which generate a scattered field $u^{sca}$.  
The aim of this work is to propose a fast and less expensive algorithm 
that can help locate all the inhomogeneous media and provide good initial guesses 
for some computationally more demanding iterative algorithms to find more accurate approximations 
of the contrast function $\chi$.  

Our algorithm will rely on the approximate contrast source obtained by backpropagation. 
Backprogation is widely used in inverse medium scatterings, see  \cite{BerBro} \cite{Levy88} 
and the references 
therein. In this section we shall give a rigorous mathematical explanation  
of the approximate contrast 
source by backpropagation. 
Let $(\cdot, \cdot)_{L^2(S)}$ and $(\cdot, \cdot)_{L^2(D)}$  be the scalar products respectively 
in $L^2(S)$ and $L^2(D)$, and $G_s^*$:\,$L^2(S)\to L^2(D)$ be the 
adjoint of operator $G_s$:\,$L^2(D)\to L^2(S)$. 
$G_s^*$ is  called {the backpropagation operator} and given by 
 \beqnx 
    (G_S^*w)(\boldsymbol x)=k^2\int_S \overline{g(\boldsymbol x,\boldsymbol x')}\,
    w(\boldsymbol x')\,ds(\boldsymbol x')
    \q \forall\,\boldsymbol x\in\mathrm{D}. 
   \eqnx 
We shall need the following backpropagation subspace of $L^2(D)$, 
$$
V_b=\mbox{span}\,\{G^*_su^{sca}\}\,,
$$ 
which is formed by all the fields generated by the backpropagation $G_s^*$  on the scattered data $u^{sca}$.
It follows from (\ref{eq:EsacMat}) that 
\begin{equation}\label{eq:EsacMat2}
    u^{sca}(\boldsymbol x)=(G_Sw)(\boldsymbol x),\quad \boldsymbol x\in\mathrm{S}\,.
\end{equation}
The backpropagation is to seek a best approximate 
solution $w_b$ to the equation (\ref{eq:EsacMat2}) in the backpropagation subspace 
$V_b$, namely 
\begin{equation}\label{eq:ls}
||u^{sca}-G_s w_b||_{L^2(S)}^2=\min_{v_b\in V_b}||u^{sca}-G_s v_b||_{L^2(S)}^2.
\end{equation}
%
It is easy to see that the solution $w_b$ to (\ref{eq:ls}) solves 
the variational system:  
\begin{equation}\label{eq:min}
 (u^{sca}-G_sw_b,G_sv_b )_{L^2(S)}=0 \quad \forall \,v_b\in V_b\,,
\end{equation}
or  equivalently, 
\begin{equation}\label{eq:min}
 (G_sw_b,G_sv_b )_{L^2(S)} =(G_s^*u^{sca}, v_b )_{L^2(D)} \quad \forall \,v_b\in V_b\,.
\end{equation}
As $w_b, v_b\in V_b$, we can write 
\begin{equation}\label{eq:w}
w_b=\lambda \,G^*_su^{sca}\,, \quad 
v_b=\mu\, G^*_su^{sca}\,
\end{equation}
for some constants $\lambda, \mu$.  Substituting the two expressions 
into  (\ref{eq:min}) we obtain 
\begin{equation}\label{eq:lam}
\lambda= \frac{||G^*_su^{sca}||^2_{L^2(D)}}{||G_sG^*_su^{sca}||^2_{L^2(S)}}\,, 
\end{equation}
which gives the approximate contrast source by backpropagation:
\begin{equation}\label{eq:wb}
w_b=\frac{||G^*_su^{sca}||^2_{L^2(D)}}
{||G_sG^*_su^{sca}||^2_{L^2(S)}}\,\,G^*_su^{sca}\,.
\end{equation}


\section{A multilevel algorithm}
\label{sec:algorithm}
    In this section we propose a fast multilevel algorithm to find the locations and geometries of 
    all the inhomogeneous media, which are described by the contrast function $\chi$ 
    in (\ref{eq:Etot}). 
    The algorithm proceeds iteratively, and carries out two important steps 
    at each iteration based on the two fundamental equations (\ref{eq:IMat}) and 
    (\ref{eq:EsacMat}), namely the state and field equations. 
    In the first step, we apply the backpropagation technique to compute 
    an approximate contrast source $w_j$ corresponding to each incident 
    $u^{inc}_j$ ($j=1,2,\cdots,N_i$). It follows from (\ref{eq:wb}) that 
    this approximation is given by
    \begin{equation}{\label{eq:Ip}}
    w_j=\frac{||G_S^* u^{sca}_j||^2_{L^2(D)}}{||G_S G_S^* u^{sca}_j||^2_{L^2(S)}}\,\,G_S^* u^{sca}_j\,,
    \quad ~~j=1, 2, \cdots, N_i\,.
    \end{equation}
    With these approximate contributions $w_j$  of the exact contrast source  $w$ corresponding 
    to each incident $u^{inc}_j$, 
    we approximate the contrast $\chi$ pointwise by minimizing the residual  
    equation corresponding to the state equation (\ref{eq:IMat}), namely 
    \begin{small}
    \begin{equation}
   \min_{\chi(\boldsymbol x)\in {\mathbf R}^1} \sum_{j=1}^{N_i}\Big| \Big(\chi u_j^{inc}-w_j+\chi G_Dw_j\Big)(\boldsymbol x)\Big|^2\,,
    \end{equation}
    \end{small}
which yields an explicit formula to compute an approximate contrast  value $\chi(\boldsymbol x)$ 
at every point $\boldsymbol x\in D$ when an approximate contrast source $w_j$ is available:
    \begin{equation}\label{eq:chi-initial}
    \chi(\boldsymbol x)
    =\frac{\sum_{j=1}^{N_i}w_j(\boldsymbol x)\,\overline{\big(u^{inc}_j+G_Dw_j)(\boldsymbol x)}}
    {\sum_{j=1}^{N_i}|(u^{inc}_j+G_Dw_j)(\boldsymbol x)|^2}\,,
    \end{equation}
    where the overbear denotes the complex conjugate.

   

Clearly  both  (\ref{eq:Ip}) and (\ref{eq:chi-initial}) are rather crude mostly, and may provide rather poor 
approximations for the exact contrast source $w$ and contrast profile $\chi$. 
But as it will be seen, 
when we combine these two poor approximations in a novel manner with some multilevel technique, 
it generates a very efficient and robust algorithm for locating  
an accurate position and shape  of each inhomogeneous medium. 

Our basic idea is also based on another simple observation. 
We know that the exact contrast function $\chi(\boldsymbol x)$  vanishes outside the scatterer $\Omega$, 
so its support provides the location and shape of 
the scatterer $\Omega$, which is formed by all the inhomogeneous media. 
This observation, along with the previous two explicit evaluation formulae (\ref{eq:Ip}) and (\ref{eq:chi-initial}) 
and a novel multilevel technique, 
forms the foundation of our new multilevel sampling algorithm. 
%

For the description of the algorithm, we first introduce two new concepts, {\it the smallest distance}
and {\it the first gap interval} with index $M$. 
For a given finite positive sequence,  $\{a_1$, $a_2$, $\cdots$, $a_m\}$, {\it its smallest distance} 
is the positive smallest one among all the distances between two neighboring elements, namely 
$dist(a_i, a_{i+1})$, $i=1, 2, \cdots, m-1$. 
Among all these $m-1$ distances, 
if there exists some $j$ such that $2\le j\le m-1$ and the distance $dist(a_j, a_{j+1})$ 
is $M$ times larger than {\it the smallest distance} 
of the sequence $\{a_1, a_2, \cdots, a_j\}$, then $[a_j, a_{j+1}]$ is called 
{\it a gap interval}. The first such interval is called {\it the first gap interval}. 


\ms
Now we are ready to state our new algorithm. 

\ss 
\textbf{\underline {Multilevel Sampling Algorithm.}}
\begin{enumerate}
\item  Choose a sampling domain $\mathbf{D}$ that contains the scatterer
$\Omega$; \\
Select a uniform (coarse) mesh on $\mathbf{D}$, consisting of rectangular (2D) or cubic  (3D) 
elements; write the mesh as $D_1$; \\
Select a tolerance $\varepsilon$ and an index $M$;
set an initial cut-off value $c_0:=0$ and $k:=1$.
\item 
    Compute an approximate value of the contrast $\chi_k(\x)$ at  each grid point 
    $\x\in D_k$, 
    using the formulae (\ref{eq:Ip}) and (\ref{eq:chi-initial}).  Then do the following:
    \begin{enumerate}
    \item[2.1]
    Order all the values  of  $\chi_k(\x)$  satisfying $\chi_k(\x)\ge c_{k-1}$ into a non-decreasing sequence; \\
    Find {\it the first gap interval} of the sequence with index $M$;  \\
    Choose the right endpoint of {\it this first gap interval} as the next cut-off value $c_k$.  
     \item[2.2]
     If $\chi_k(\x)\ge c_k$ at a grid point ${\x}$, select all the grid 
     points of the elements which share ${\x}$ as one of its vertices; \\
     Remove all the grid points in $D_k$, which are not selected; \\
     Update $D_k$ by all those selected grid points. 
     

    \end{enumerate}
\item If $|c_k-c_{k-1}|\le \varepsilon$, go to Step 4; \\
otherwise refine the mesh $D_k$ to get $D_{k+1}$; set $k:=k+1$ and go to Step 2. 

\item Output all grid points in $D_k$ for 
 the supports of all inhomogeneous media; \\
  Output the approximate contrast values $\chi_k(\x)$ at all grid points of $D_k$. 
\end{enumerate}

%
We can easily see that the Multilevel Sampling Algorithm does not involve any optimization process 
or matrix inversions, and its major cost is to update the contrast values 
using the explicit formulae (\ref{eq:Ip}) and (\ref{eq:chi-initial})
at each iteration, and the computational sampling domain $D_k$ shrinks 
as the iteration goes. So the algorithm is rather simple and less expensive. 
In addition, as the cut-off values are basically to distinguish the homogeneous background medium 
where $\chi(\x)$ vanishes and  the inhomogeneous media where $\chi(\x)$ should be essentially 
different from $0$, so our cut-off values are rather easy to choose 
and insensitive to the size and physical features 
of scatterers. In fact, the cut-off value can start simply with zero, then 
it is updated automatically with the iteration. As we shall see from numerical examples 
in the next section, 
the algorithm works well with few incidents, even with one; and it is self-adaptive, namely 
it can recover some elements that have been removed 
at the previous iterations due to the computational errors. 
In terms of these aspects the new Multilevel Sampling Algorithm outperforms 
the popular linear sampling methods \cite{Pot}, including 
the improved multilevel variant \cite{LLZ08}.

{\bf Discretization}. We end this section with a brief discussion about 
numerical discretization of the integrations involved in the above 
multilevel algorithm.  
   We illustrate only the discretization of   (\ref{eq:IMat}) and (\ref{eq:EsacMat}), 
   as all other integrations involved in the algorithm can be approximated similarly.  
   To do so, we divide the domain $\mathbf{D}$ into smaller rectangular or cubic 
   elements, whose centers 
   are denoted as $\boldsymbol x_1,\boldsymbol x_2,\dots,\boldsymbol x_M$. 
   Using the coupled-dipole method (CDM) or discrete dipole approximation (DDA) \cite{BelChaSen,Lak},
   we can discretize (\ref{eq:IMat}) by 
    \begin{equation}\label{eq:EtotExp2}
    w_j(\boldsymbol x_m)=\chi(\boldsymbol x_m) u^{inc}_j(\boldsymbol x_m)
    +k^2\chi(\boldsymbol x_m)\sum_{n\neq m}A_{n}g(\boldsymbol x_m,\boldsymbol x_n)w_j(\boldsymbol x_n), \quad m=1,2,\dots,M,                                                
    \end{equation}
    where $A_{n}$ is the area or volume of the $n$-th element. 
    Similarly, we can discretize equation (\ref{eq:EsacMat}) at every point 
    $\boldsymbol x\in\mathrm{S}$ by 
    \begin{eqnarray}\label{eq:Esca2}
    u^{sca}_j(\boldsymbol x)=k^2\sum_{n=1}^M A_ng(\boldsymbol x,\boldsymbol x_n)
    w_j(\boldsymbol x_n) \q \mbox{for} \q j=1,2,\cdots,N_i.
    \end{eqnarray}

\section{Numerical simulations}
In this section we present several examples to test the effectiveness of 
our multilevel sampling algorithm.
     We first list the parameters that are used in our  numerical simulations: 
    the wave number $k$ and wave length $\lambda$ are taken to be $k=2\pi$ and $\lambda=1$; 
    the number of incidences and receivers are set to be  $N_i=6$ and $N_s=30$ respectively, 
    both equally distributed on the circle of radius $5\lambda$; 
     the index $M$ of {the first gap interval} 
     and the tolerance parameter $\varepsilon$ are chosen to be $100$ and $10^{-3}$ respectively. 
     In all the numerical simulations, random noises are added to the exact scattering data
     in the following form:
         \[u^{sca}_j(\boldsymbol x): = u^{sca}_j(\boldsymbol x)[1+\xi(r_{1,j}(\boldsymbol x)+i r_{2,j}(\boldsymbol x))], \quad j=1,2,\cdots,N_i\]
     where $r_{1,j}(\boldsymbol x)$ and $r_{2,j}(\boldsymbol x)$ are two random numbers varying 
     between -1 and 1, and $\xi$ corresponds to the level of the noise, which is usually 
     taken to be $10\%$ unless specified otherwise. 
     All the programs in our experiments are written in MATLAB and run on a 2.83GHz PC with 4GB memory.
     
\subsection{Two-dimensional reconstructions}
\label{sec:numerics}

\begin{figure}[!htb]
 \hfill{}\includegraphics[clip,width=0.45\textwidth]{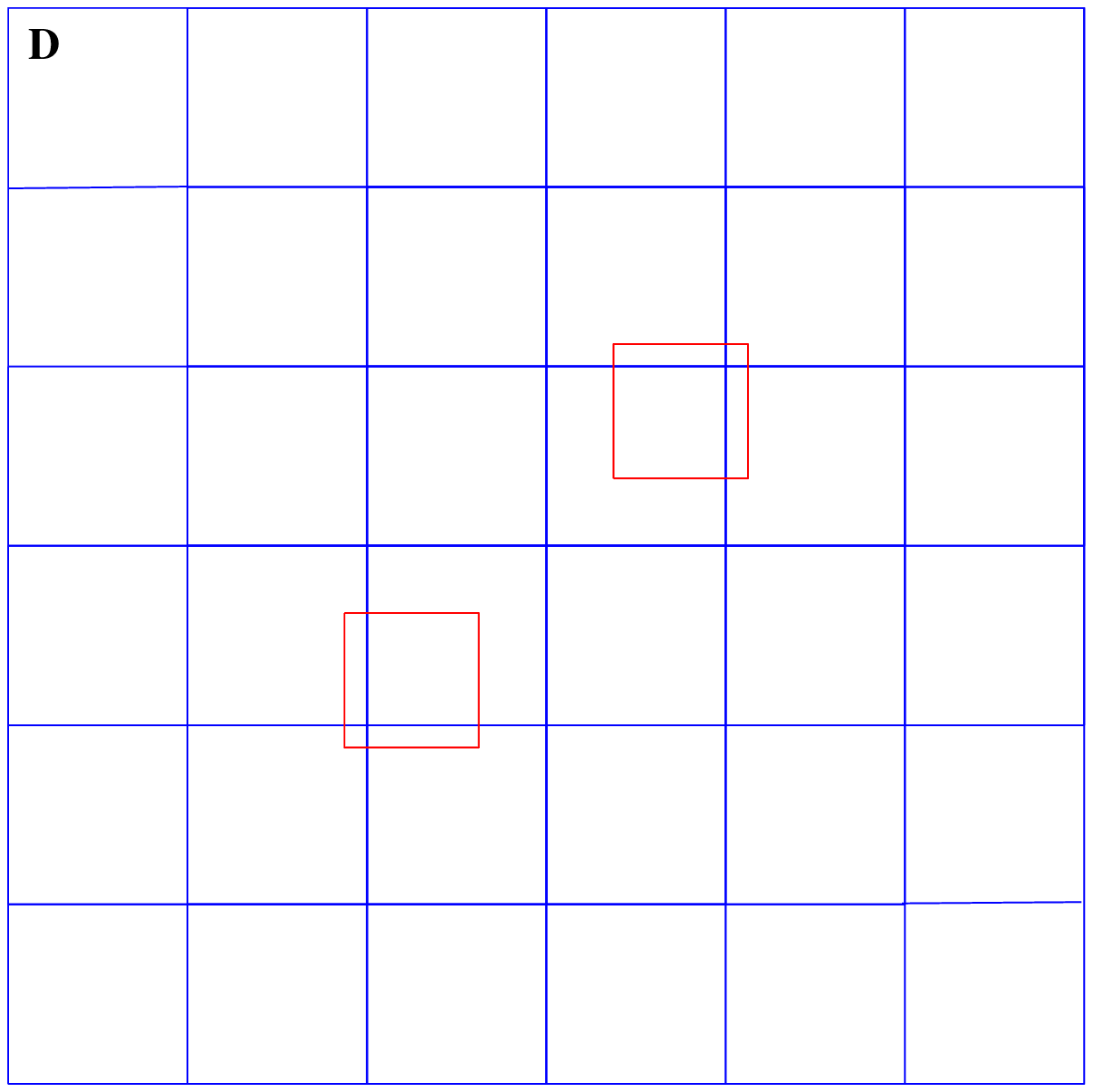}\hfill{}
 \hfill{}\includegraphics[clip,width=0.45\textwidth]{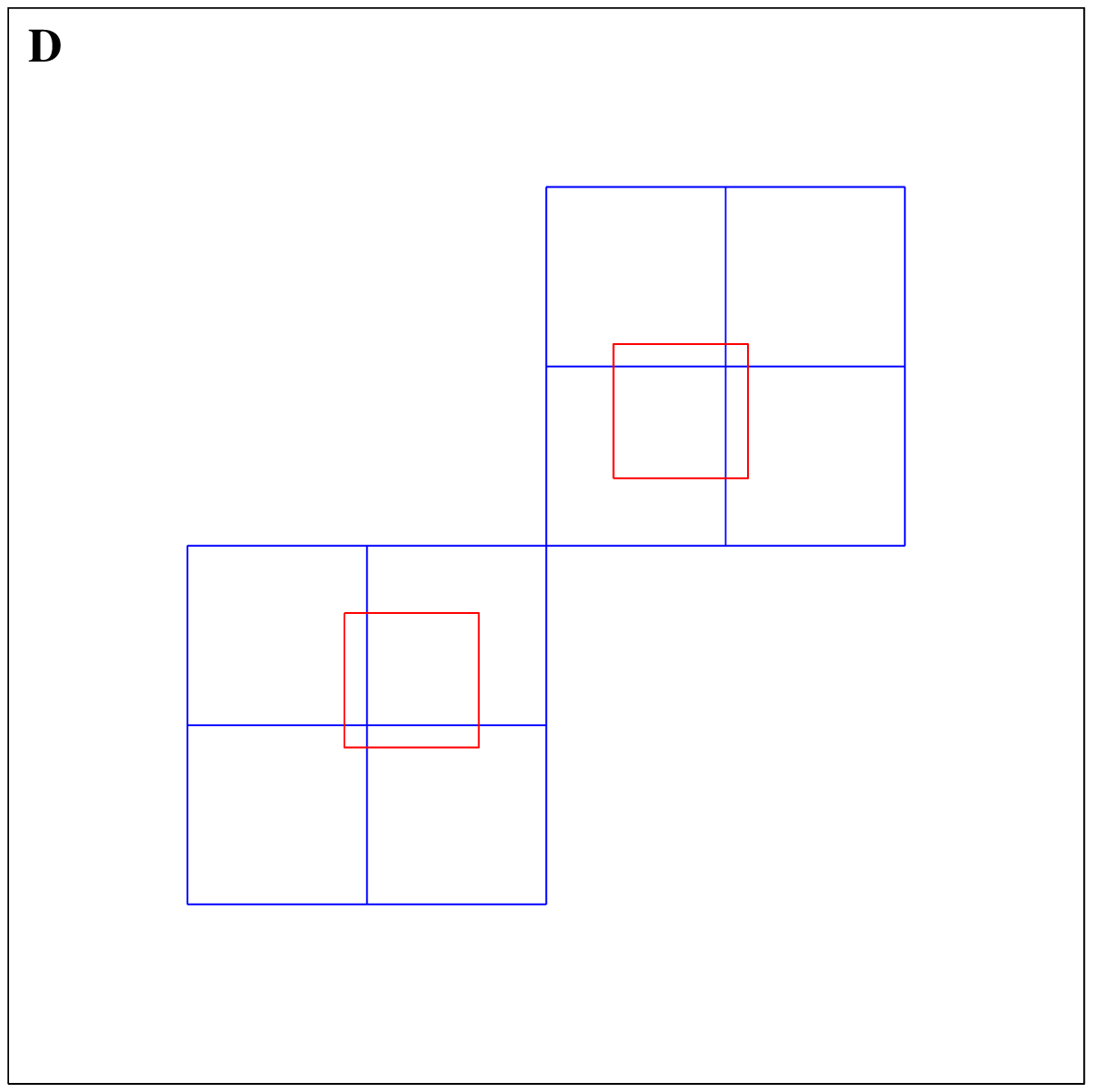}\hfill{}
 
 \hfill{}(a)\hfill{} \hfill{}(b)\hfill{}
 
 \hfill{}\includegraphics[clip,width=0.45\textwidth]{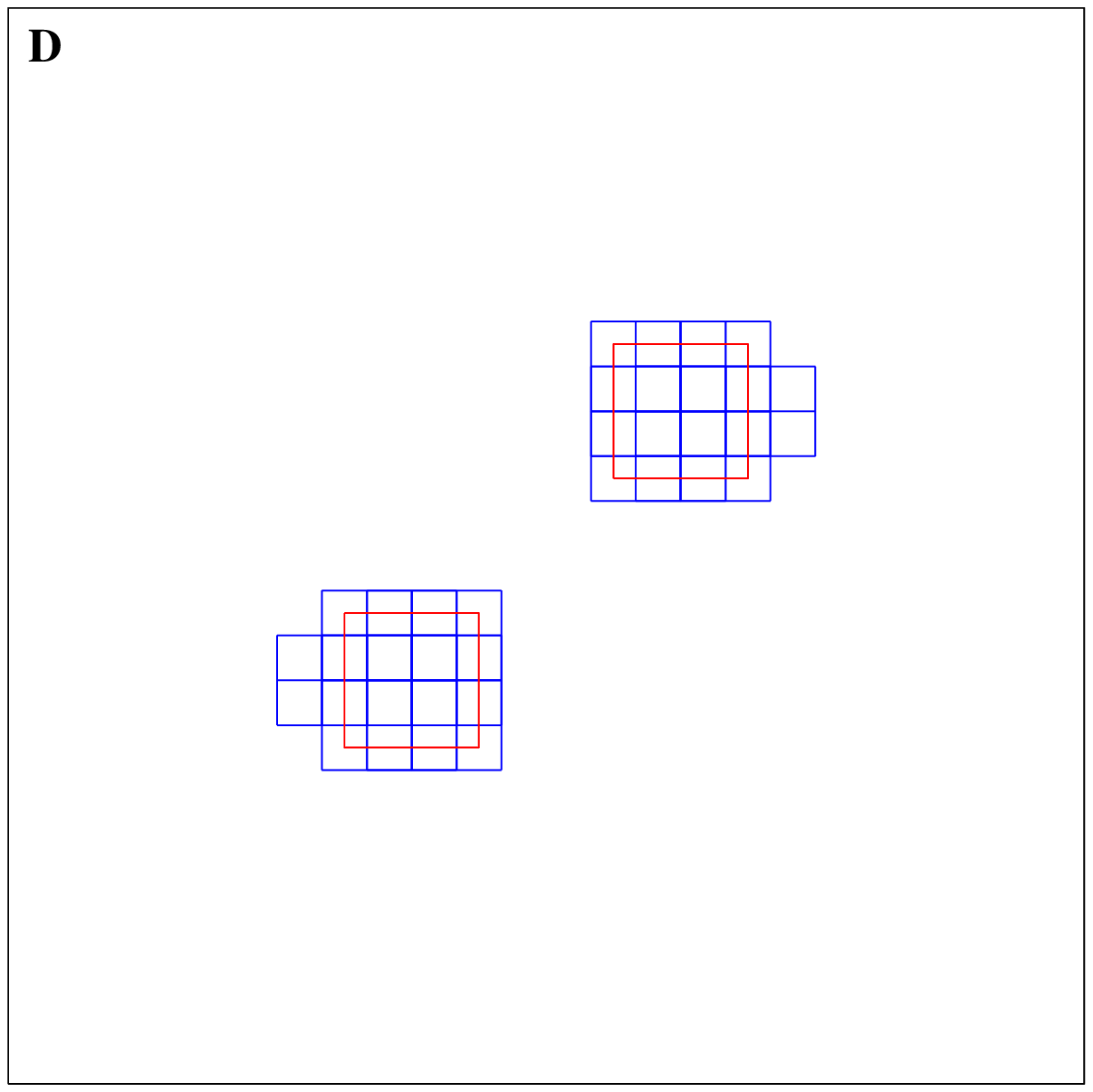}\hfill{}
 \hfill{}\includegraphics[clip,width=0.45\textwidth]{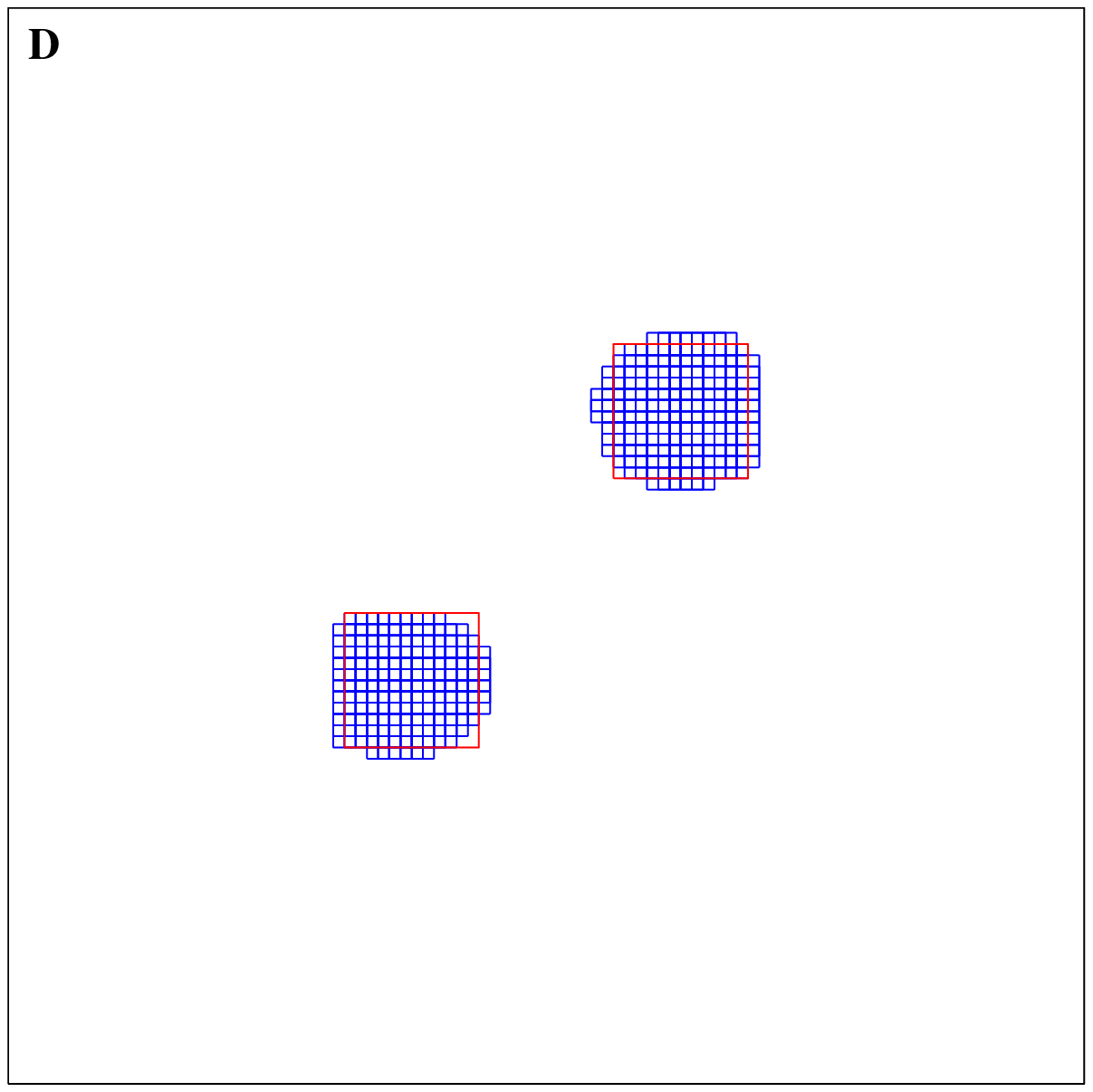}\hfill{}

\hfill{}(c)\hfill{} \hfill{}(d)\hfill{}

%
%
 \caption{\label{fig:simulation1}\emph{(a) The initial (coarse) mesh on the sampling domain; 
 (b)-(d) Reconstructions at the 1st, 3rd and 5th iteration.}}
 \end{figure}

{\bf Example 1}. This example shows a scatterer $\Omega$ consisting of two squares of side length 
$0.3\lambda$, located respectively at $(-0.3\lambda,-0.3\lambda)$
 and $(0.3\lambda,0.3\lambda)$, with 
 their contrast values being 1 and 2 respectively; see the two red squares 
 in Figure~\ref{fig:simulation1}(a).  We take the sampling domain 
 $\mathbf{D}=[-1.2\lambda,1.2\lambda]\times[-1.2\lambda,1.2\lambda]$, which is quite large  
 compared to the scatterer $\Omega$, with its area of 64 times of the area of one scatterer  
 component.   More importantly, we see 
 that these two small objects are quite close to each other. 
The mesh refinement during the multilevel algorithm is  carried out 
based on the rule
$h_{k}={0.4\lambda}/{2^{k}}$, where $k$ is the $k$-th refinement, and $h_0$ and $h_k$ 
are respectively the mesh sizes of the initial mesh and the mesh after the $k$th refinement.  
 The numerical reconstructions are shown 
 in Figure~\ref{fig:simulation1}(b)~-~\ref{fig:simulation1}(d) respectively for the 1st, 3rd and 5th iteration. One can observe from the figures that the algorithm converges 
 very fast 
 and provides accurate locations of the two medium components in only 5 iterations. 
 Moreover, we can see an important advantage of the algorithm, 
 i.e., it can separate the disjoint medium components  quickly.
  
\begin{figure}[!htb]
  \hfill{}\includegraphics[clip,width=0.45\textwidth]{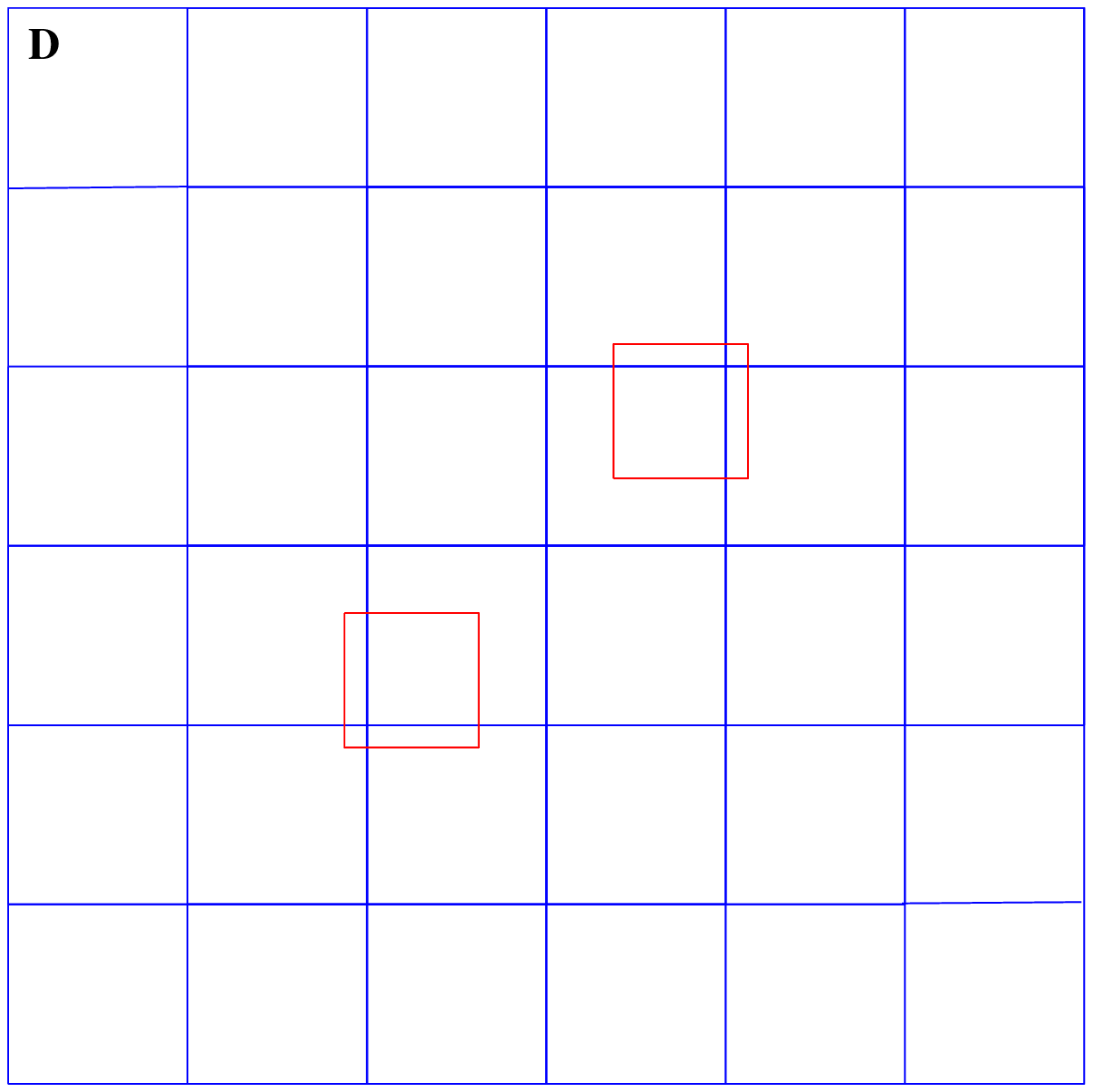}\hfill{}
 \hfill{}\includegraphics[clip,width=0.45\textwidth]{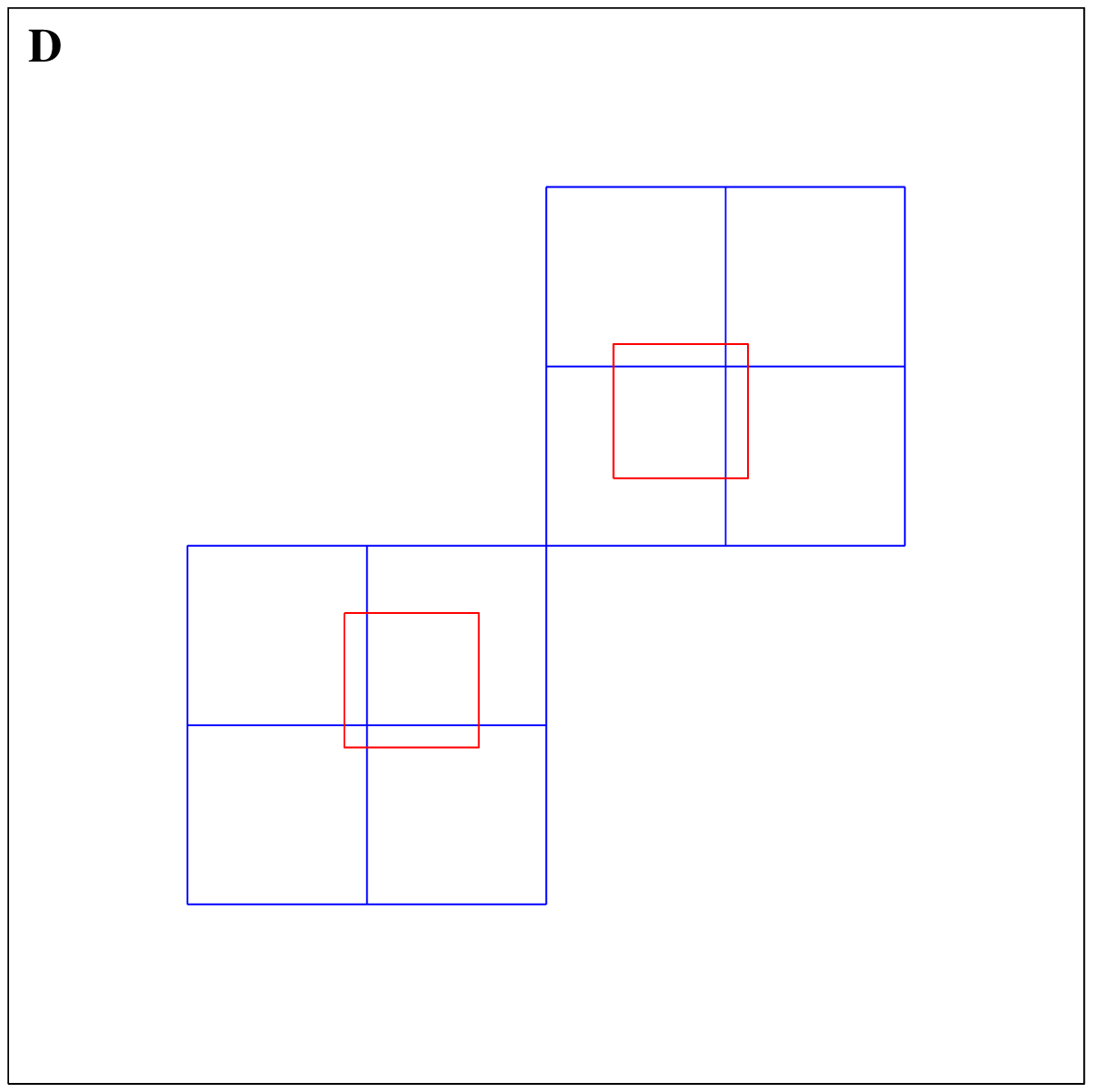}\hfill{}

 \hfill{}(a)\hfill{} \hfill{}(b)\hfill{}
 
 \hfill{}\includegraphics[clip,width=0.45\textwidth]{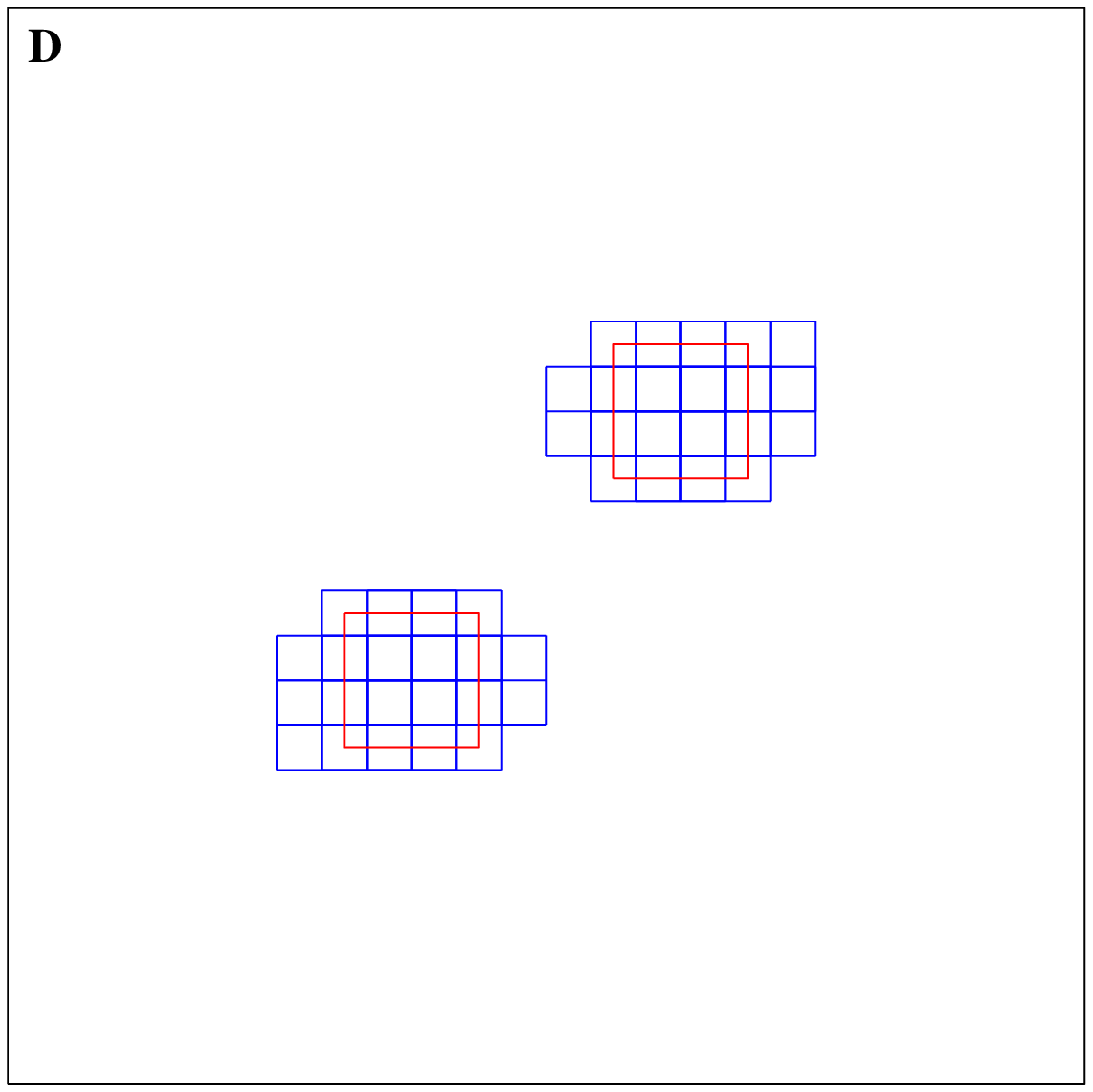}\hfill{}
 \hfill{}\includegraphics[clip,width=0.45\textwidth]{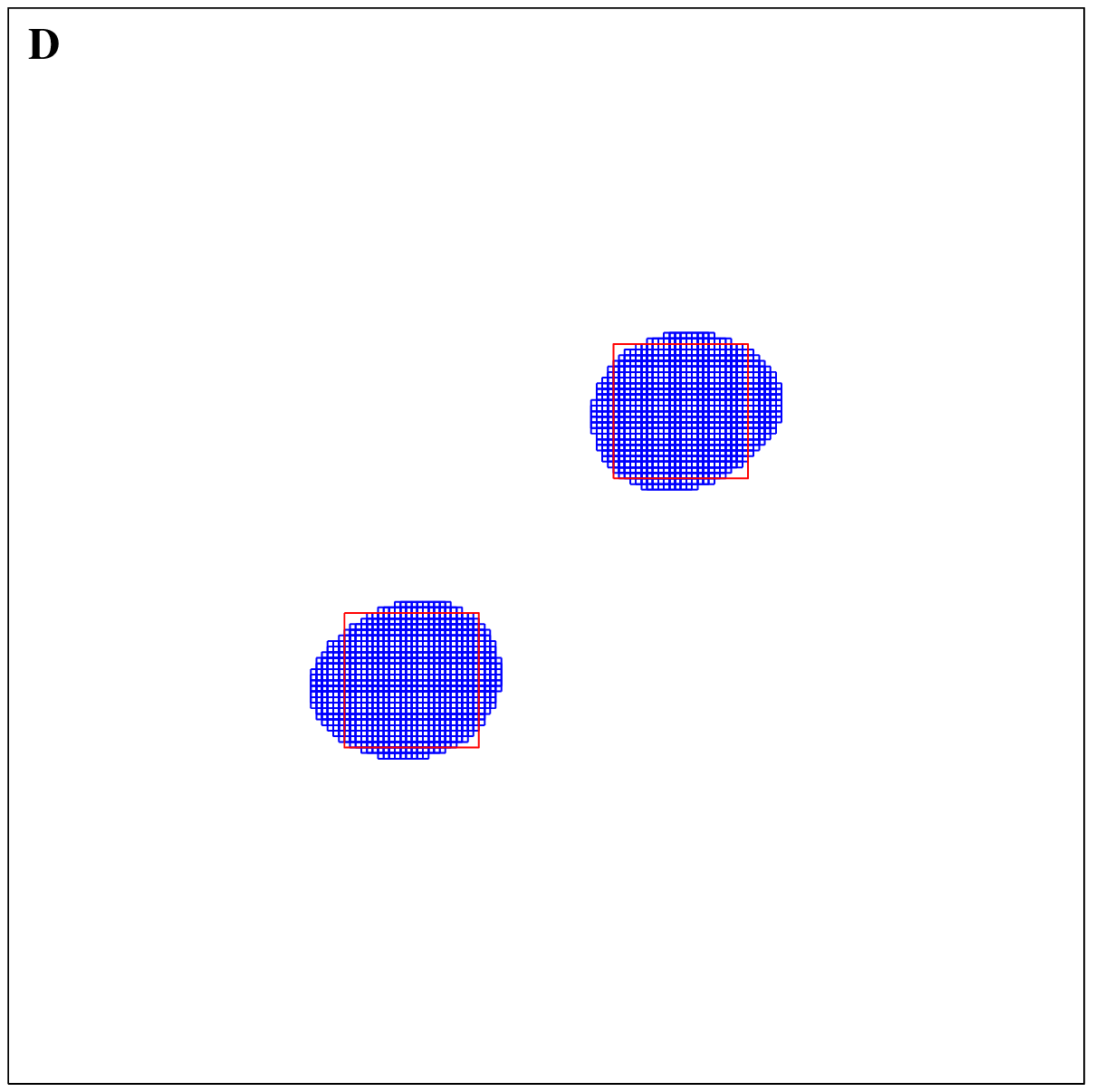}\hfill{}

 \hfill{}(c)\hfill{} \hfill{}(d)\hfill{}

%

 \caption{\label{fig:simulation2}\emph{(a) The first coarse mesh on the sampling region; (b)-(d) Reconstructions at the 1st, 3rd and 6th iteration.}}
 \end{figure}

\ss
{\bf Example 2}. This example is the same as Example 1, 
except that the contrast values of two medium components are now variable functions, namely 
$$
\chi(x, y)=sin\frac{\pi(10|x|-1.5)}{3}sin\frac{\pi(10|y|-1.5)}{3}.
$$
The numerical reconstructions are shown in Figure~\ref{fig:simulation2}(b)-~\ref{fig:simulation2}(d)
for the 1st, 3rd and 6th iteration. Again, we observe from the figures that the algorithm converges fast,
 provides very satisfactory locations of the two medium components in only 6 iterations, 
 and it can separate the disjoint medium components  quickly.

\begin{figure}[!htb]
 \hfill{}\includegraphics[clip,width=0.45\textwidth]{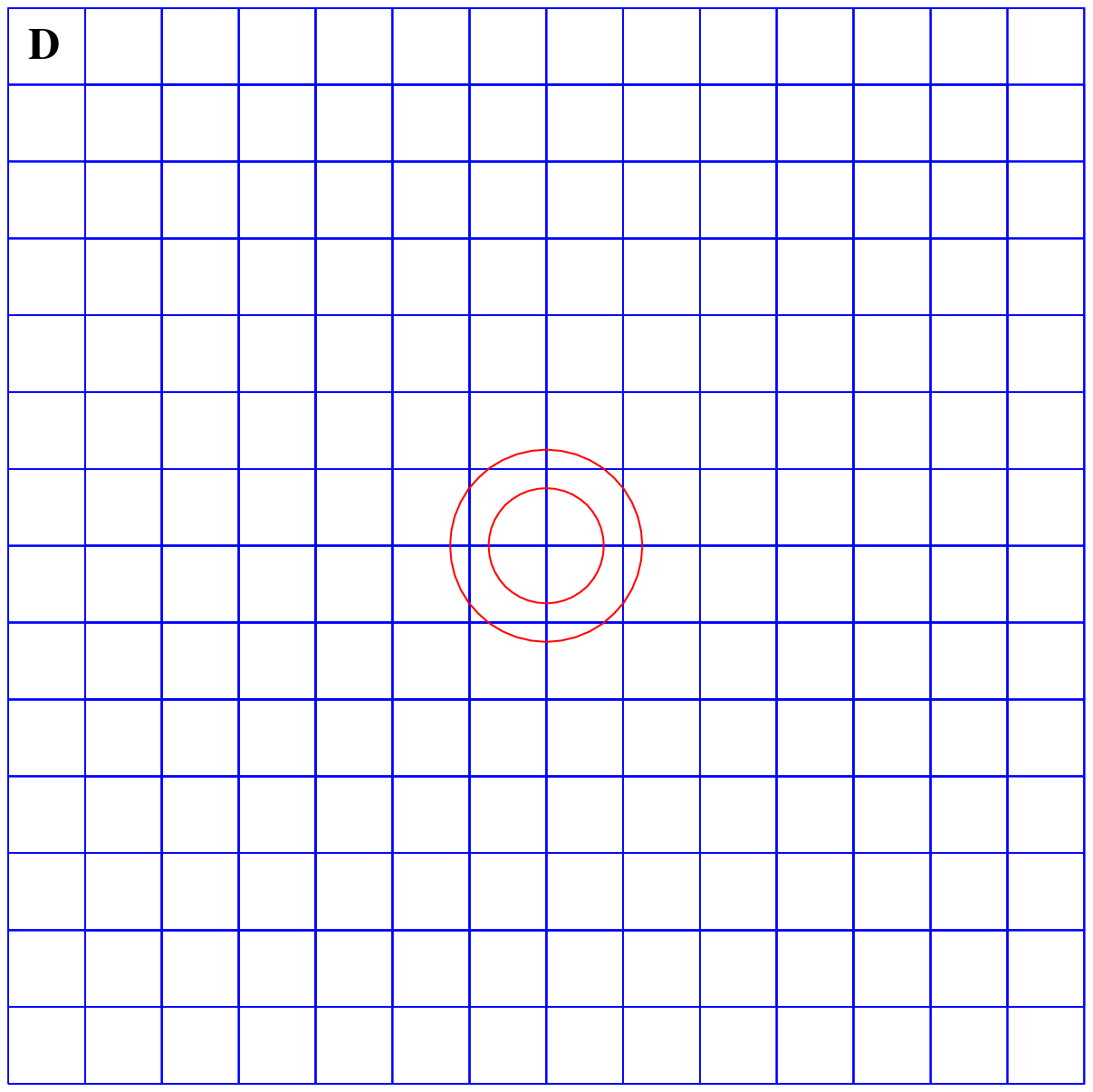}\hfill{}
 \hfill{}\includegraphics[clip,width=0.45\textwidth]{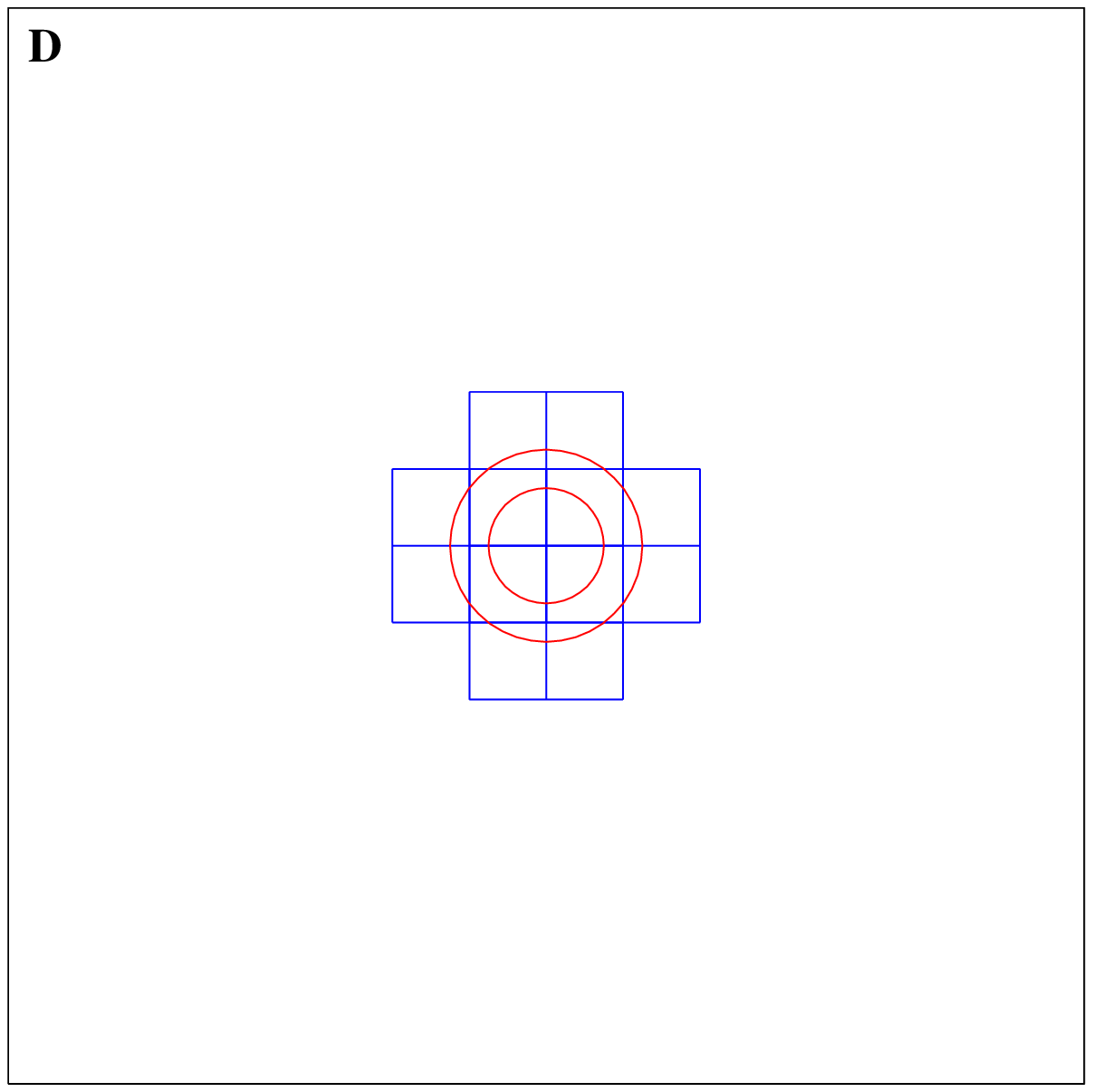}\hfill{}

 \hfill{}(a)\hfill{} \hfill{}(b)\hfill{}
 
 \hfill{}\includegraphics[clip,width=0.45\textwidth]{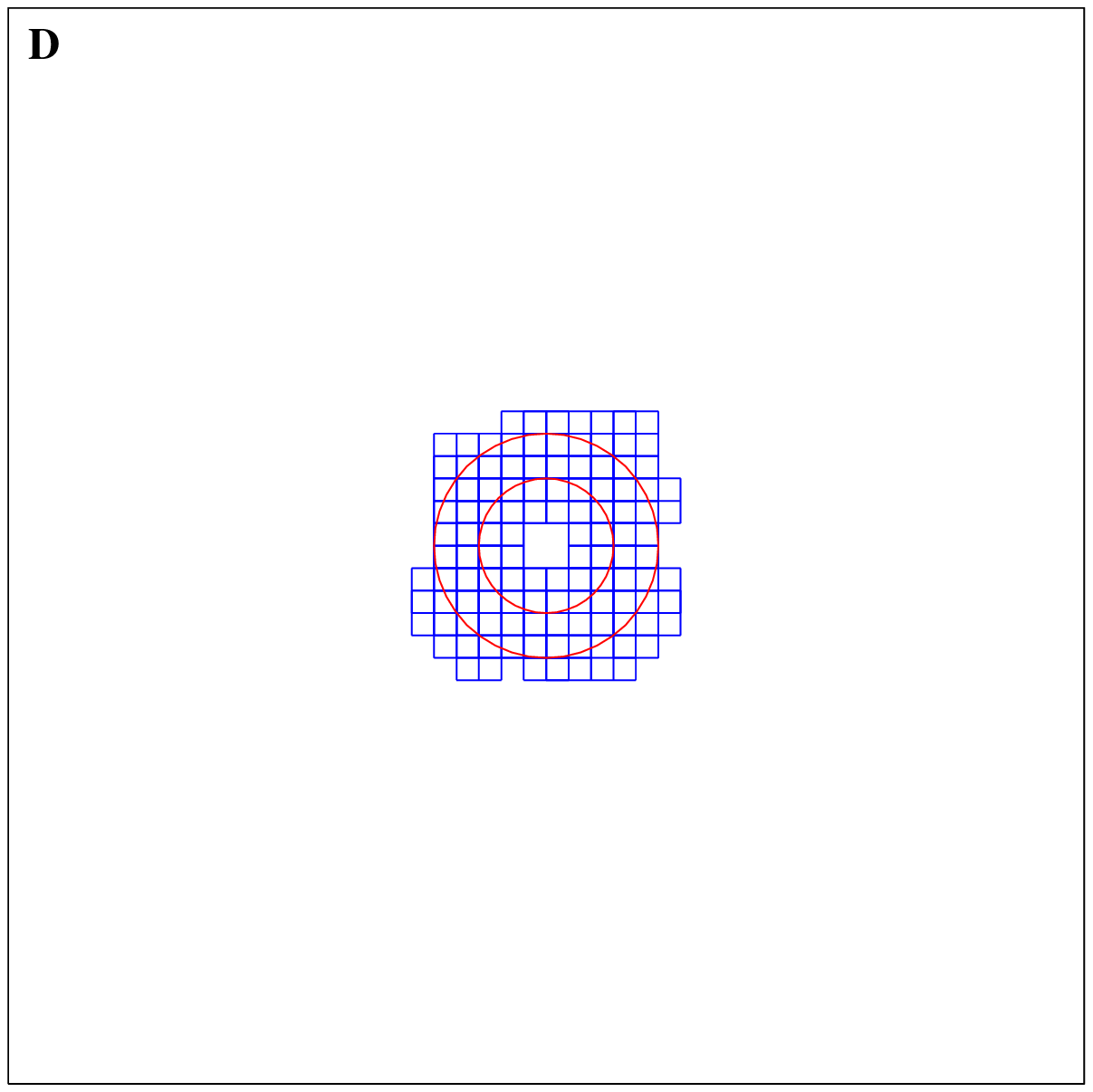}\hfill{}
 \hfill{}\includegraphics[clip,width=0.45\textwidth]{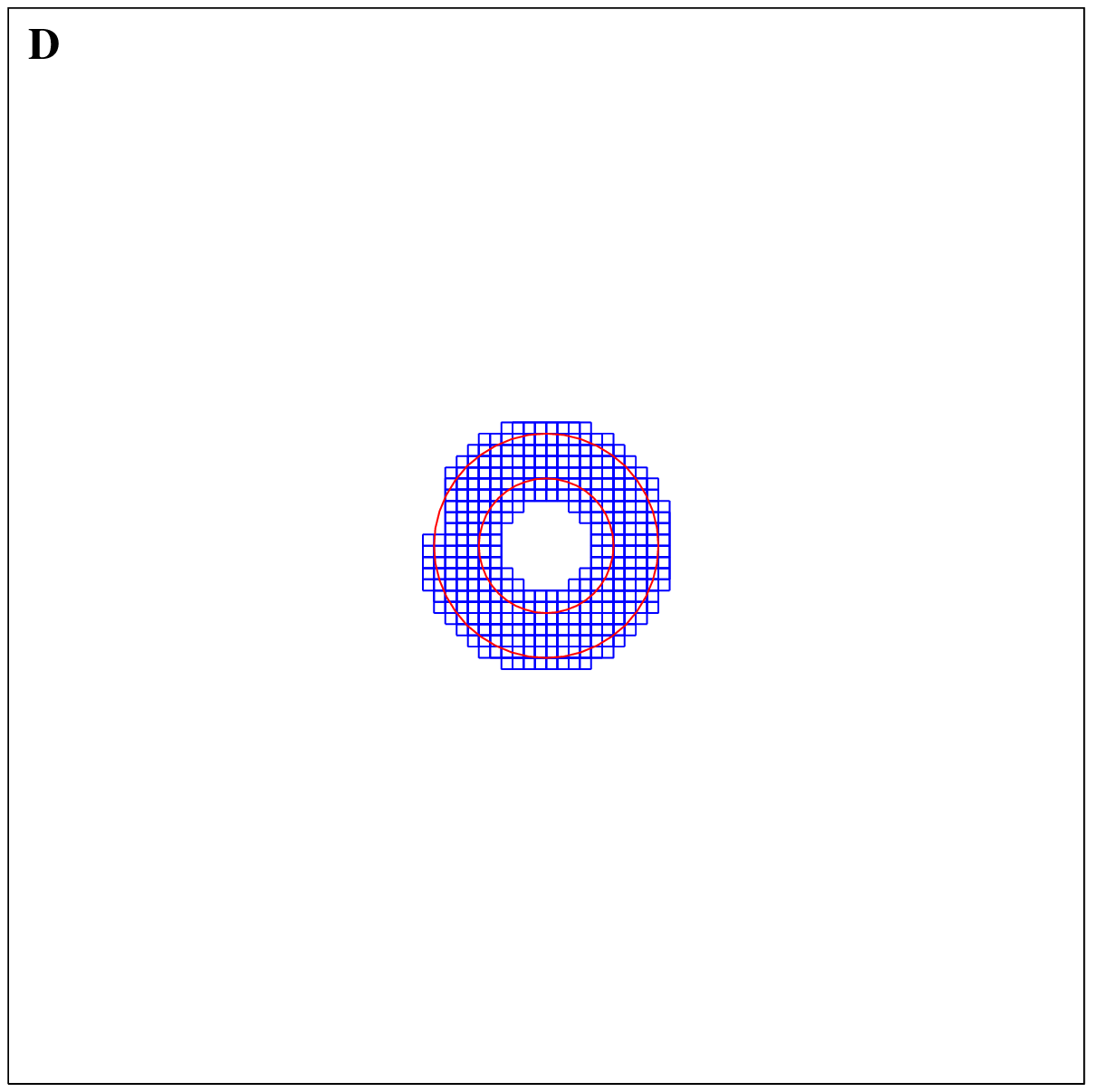}\hfill{}

 \hfill{}(c)\hfill{} \hfill{}(d)\hfill{}

%

 \caption{\label{fig:simulation3}\emph{ (a) The first coarse mesh on the sampling region; (b)-(d) 
 Reconstructions at the 1st, 3rd and 4th iteration.}}
 \end{figure}

\ss
{\bf Example 3}. This example considers a scatterer $\Omega$ of a thin annulus 
with the inner and outer radii being $0.3\lambda$ and $0.5\lambda$ respectively 
and centered at the origin. The contrast value $\chi(\boldsymbol x)$ is $2$ inside the  thin annulus.
The sampling domain $\mathbf{D}$  is taken to be a square of side length with $5.6\lambda$, 
as shown in Figure~\ref{fig:simulation3}(a).

It is easy to see the sampling domain $\mathbf{D}$ has an area 
of  about 62 times as large as the one of the annulus, 
and the annulus has a very thin thickness, i.e., $0.2\lambda$. 
The mesh refinement during the multilevel algorithm is  carried out 
based on the rule
$h_{k}={0.4\lambda}/{2^{k}}$, where $k$ indicates the $k$-th refinement and $h_k$ 
is the mesh size after the $k$th refinement.  
 The numerical reconstructions are shown in Figure~\ref{fig:simulation3}(b)-~\ref{fig:simulation3}(d)
 for the 1st, 3rd and 4th iteration.
 Same as for the previous two examples, 
 the reconstructions are quite satisfactory and the accurate locations of 
 the scatterer can be achieved.    
   
\subsection{Reconstruction for the contrast function $\chi$}

   There are many numerical reconstruction methods available in the literature 
   for the contrast profile function $\chi$,  that are 
   robust and more accurate than our new multilevel method. But these methods 
   are usually more complicated technically and much more demanding computationally, 
   as they mostly involve nonlinear optimizations and matrix inversions. Without a reasonably good 
   initial location for each inhomogeneous medium, we may have to take a much larger 
   sampling domain than required so these methods can be extremely time 
   consuming, especially in three dimensions. 
   Using the newly proposed multilevel algorithm in Section \ref{sec:algorithm}, we can 
   first locate a much smaller sampling domain than usual (or the one we originally selected) 
   in a numerical reconstruction for the contrast $\chi$. Then we can apply any existing reconstruction 
   algorithms for more accurate reconstructions, starting with an initial sampling domain provided 
   by the new multilevel algorithm. This may save us a great fraction of the entire computational 
   costs. 
   
\begin{figure}[!htb]
 \hfill{}\includegraphics[clip,width=0.45\textwidth]{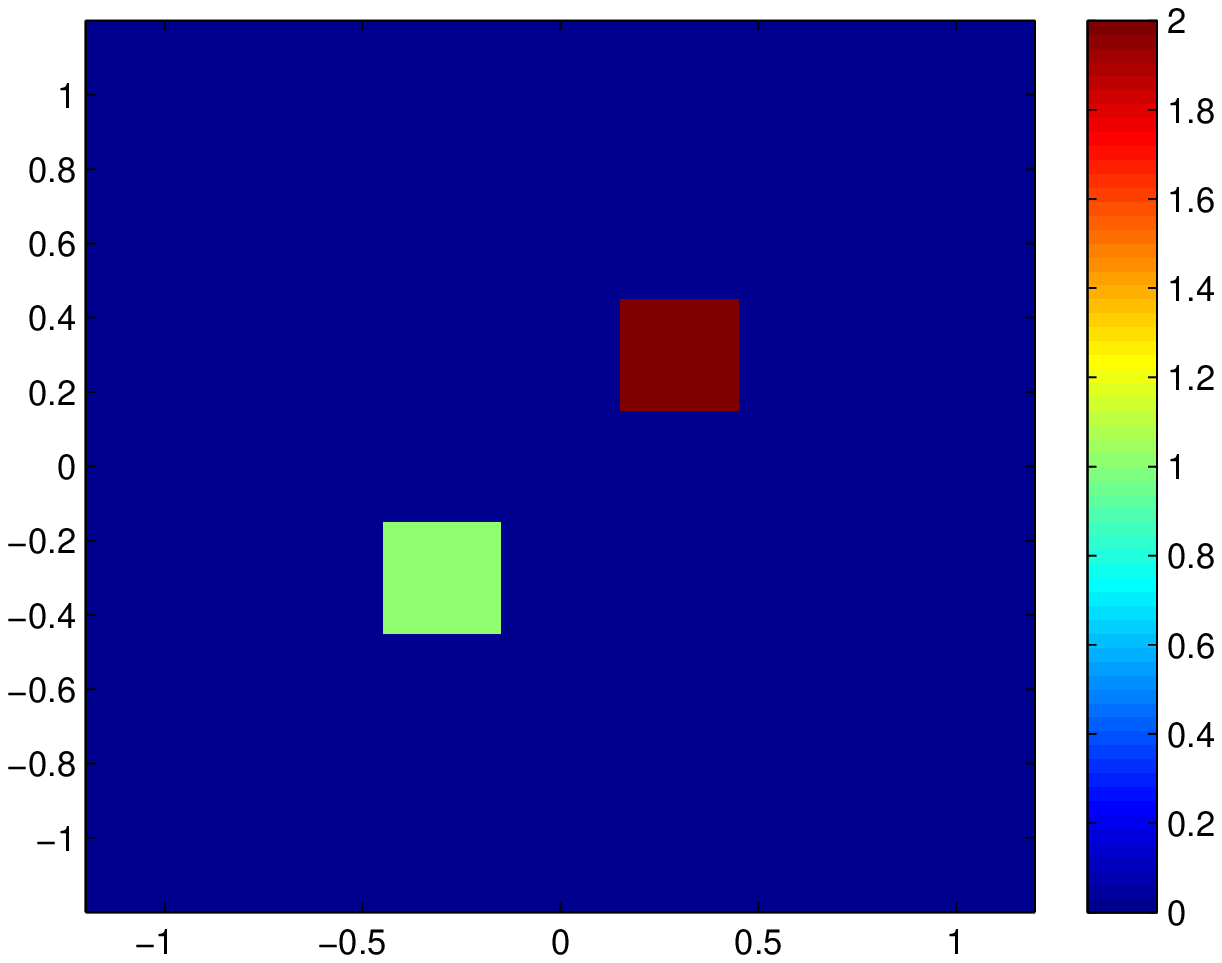}\hfill{}
 \hfill{}\includegraphics[clip,width=0.45\textwidth]{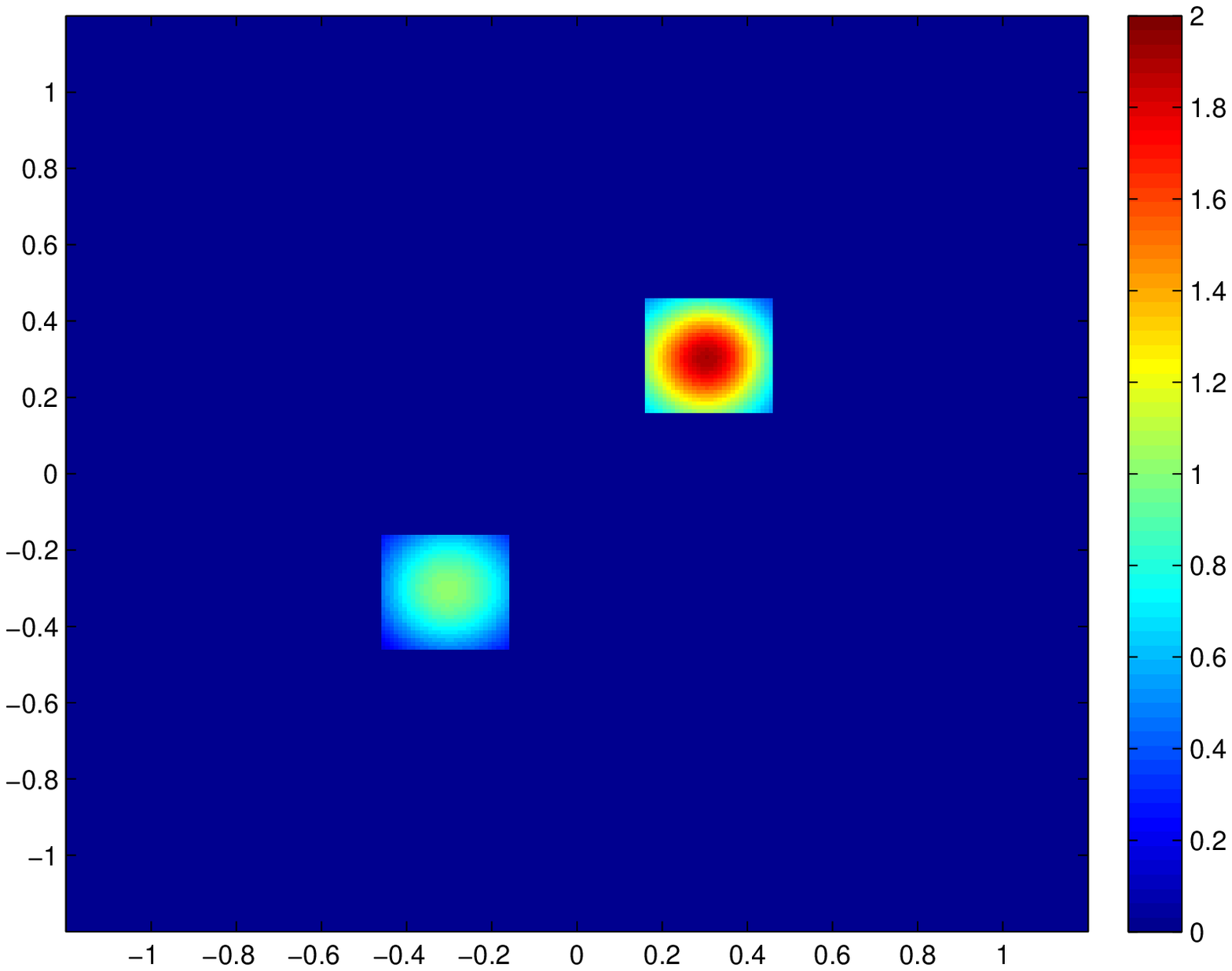}\hfill{}

 \hfill{}(a)\hfill{} \hfill{}(b)\hfill{}
 
 \hfill{}\includegraphics[clip,width=0.45\textwidth]{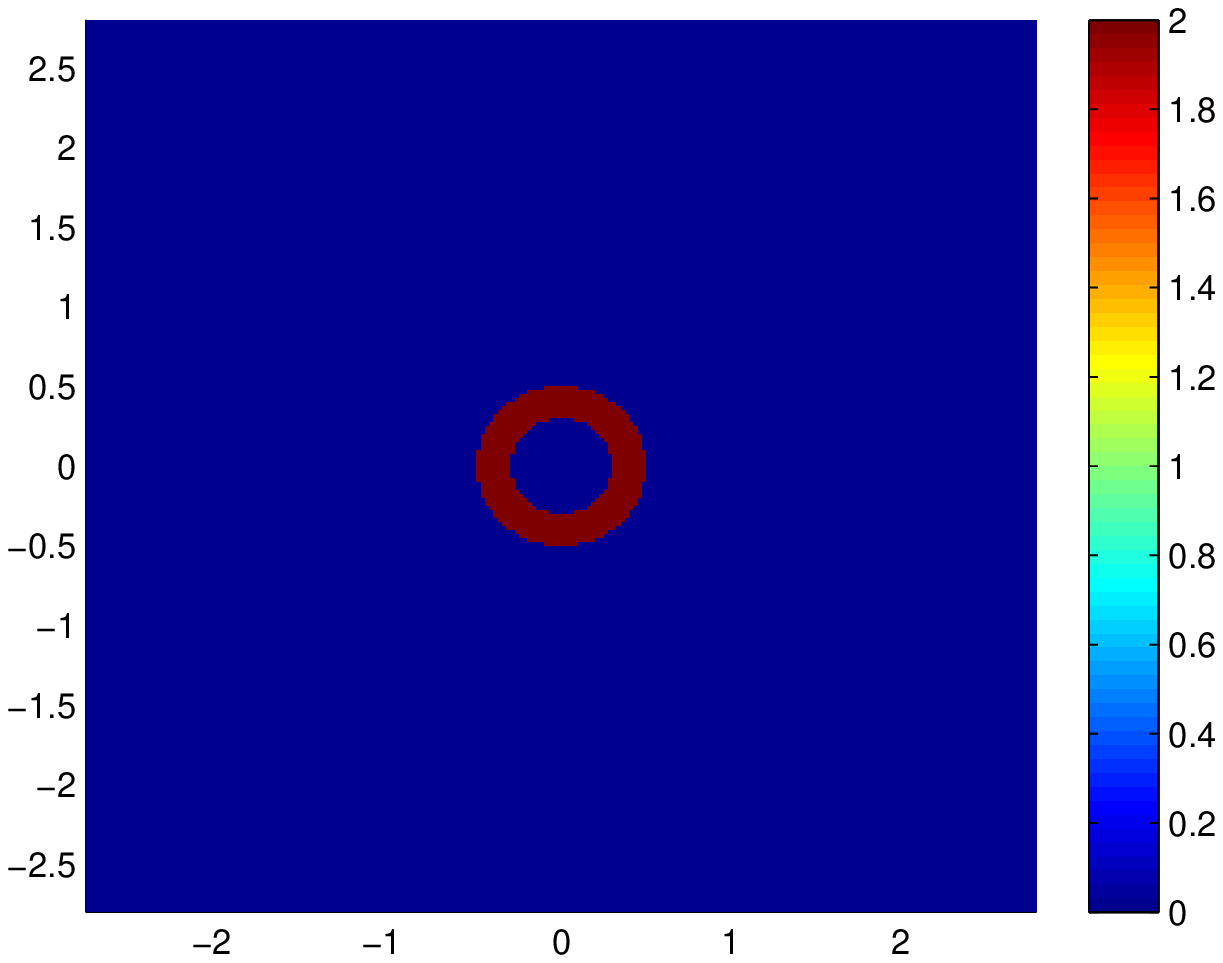}\hfill{}
 \hfill{}\includegraphics[clip,width=0.45\textwidth]{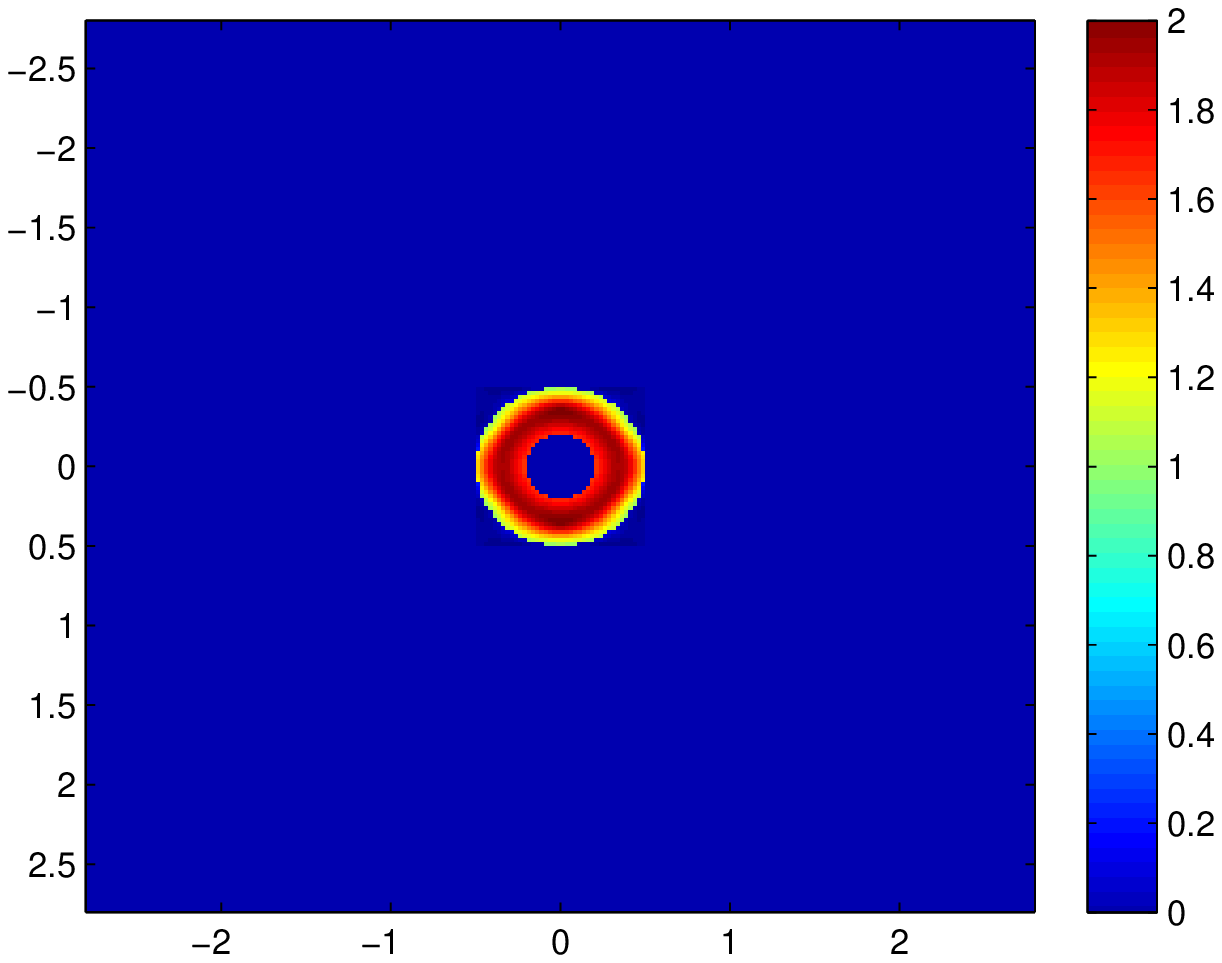}\hfill{}

 \hfill{}(c)\hfill{} \hfill{}(d)\hfill{}

 \caption{\label{fig:simulation}\emph{(a) True scatterer in Example 1;
  (b) Reconstruction by combining the multilevel method with ECSI for Example 1; 
  (c) True scatterer in Example 3; (d) Reconstruction by combining the multilevel method 
  with ECSI for Example 3}}
 \end{figure}

   Next we show numerical tests by combining our new multilevel method with 
   the popular extended contrast source inversion (ECSI) method \cite{BerBro}. 
   Firstly, we consider the same scatterer $\Omega$ and the set-ups as in Example 1 of Section \ref{sec:numerics}; see Figure \ref{fig:simulation}(a). 
   Then we apply the ECSI method \cite{BerBro} with mesh size $h=0.01\lambda$ 
   to the reconstructed domain (cf.~Figure \ref{fig:simulation1}(d)) by the multilevel algorithm. The reconstruction result is shown in Figure \ref{fig:simulation}(b). Clearly, one can see that both the location and the values of the contrast are well reconstructed.

   Secondly, we consider the same scatterer $\Omega$ and the set-ups as in Example 3 of Section \ref{sec:numerics}; see Figure (\ref{fig:simulation})(c). 
   We apply the ECSI method \cite{BerBro} with mesh size $h=0.025\lambda$ to the reconstructed domain 
   (cf.~Figure \ref{fig:simulation3}(d)) obtained by the multilevel algorithm. 
   The reconstruction result is shown in Figure \ref{fig:simulation}(d). Clearly, 
   both of the location and the values of the contrast are well reconstructed.
 
\subsection{Three-dimensional reconstructions}

\begin{figure}[!htb]
 \hfill{}\includegraphics[clip,width=0.45\textwidth]{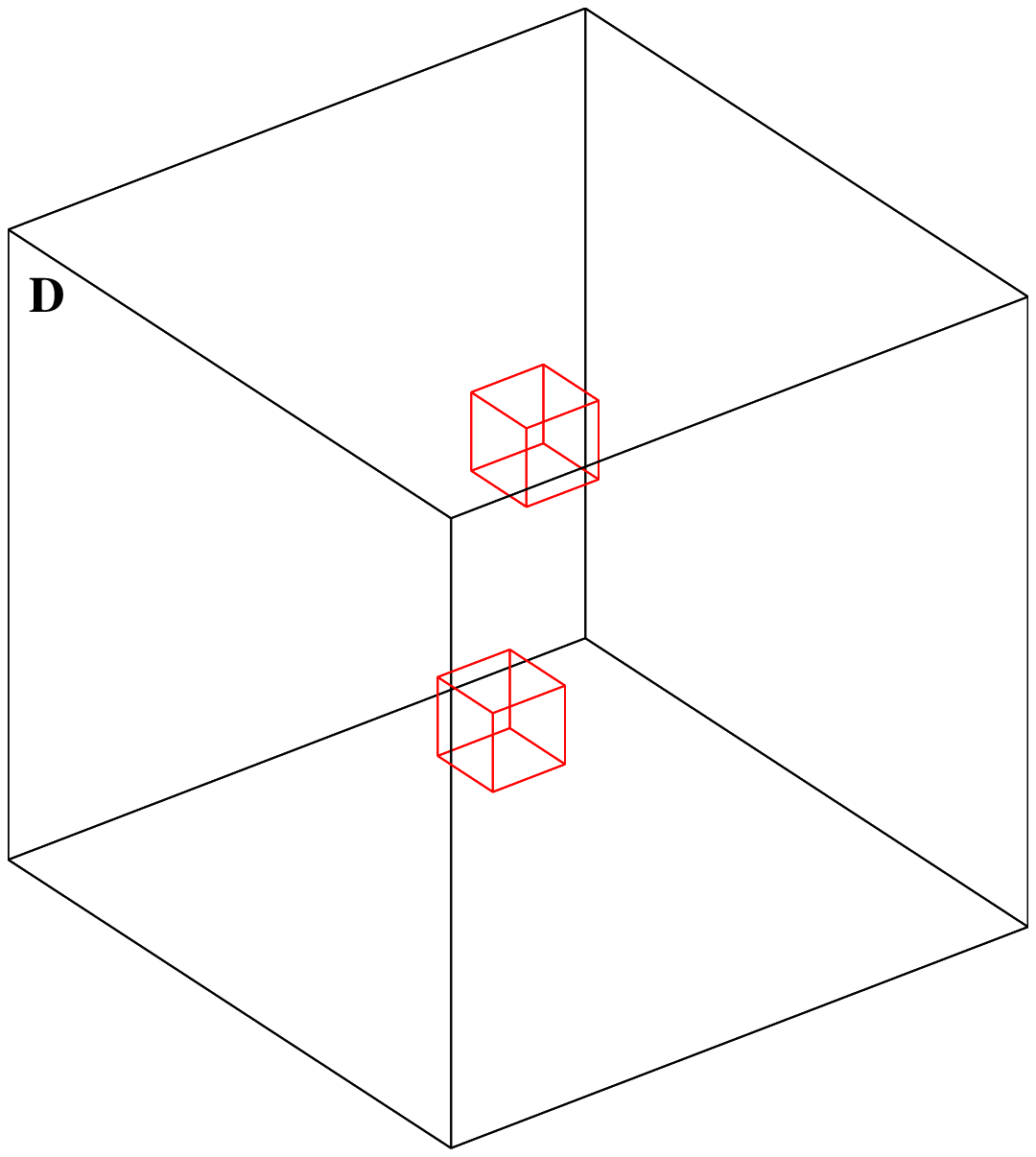}\hfill{}
 \hfill{}\includegraphics[clip,width=0.45\textwidth]{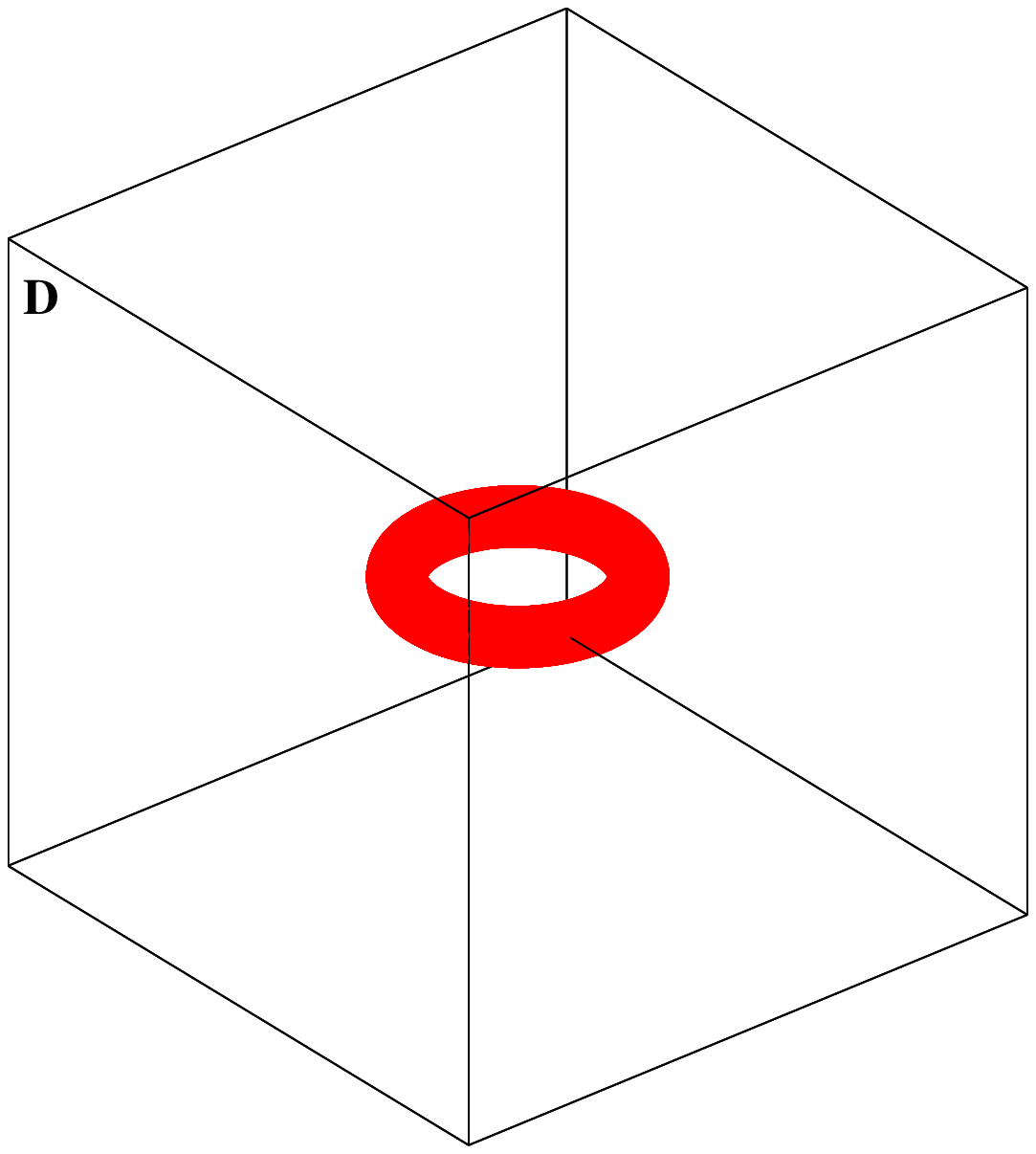}\hfill{}

 \hfill{}(a)\hfill{} \hfill{}(b)\hfill{}

 \caption{\label{fig:boxes}Scatterers  imbedded in a large sampling domain: 
 two cubic components close to each other (a); 
 a torus (b)}
 \end{figure}

{\bf Example 4}. This example tests a three-dimensional scatterer 
$\Omega$ consisting of two small cubic components:
$$
\Omega_1=[-0.45\lambda, ~-0.15\lambda]^3\,, 
\q 
\Omega_2 = [0.15\lambda, ~0.45\lambda]^3\,.
$$
The two squares are quite close to each other, both with constant contrast values $2$;
see Figure~\ref{fig:boxes}(a). 
We take the sampling domain to be $\mathbf{D}=[-1.2\lambda,1.2\lambda]^3$, which is 
about 500 times of the volume of $\Omega_1$ or $\Omega_2$. 

\begin{figure}[!htb]
 \hfill{}\includegraphics[clip,width=0.45\textwidth]{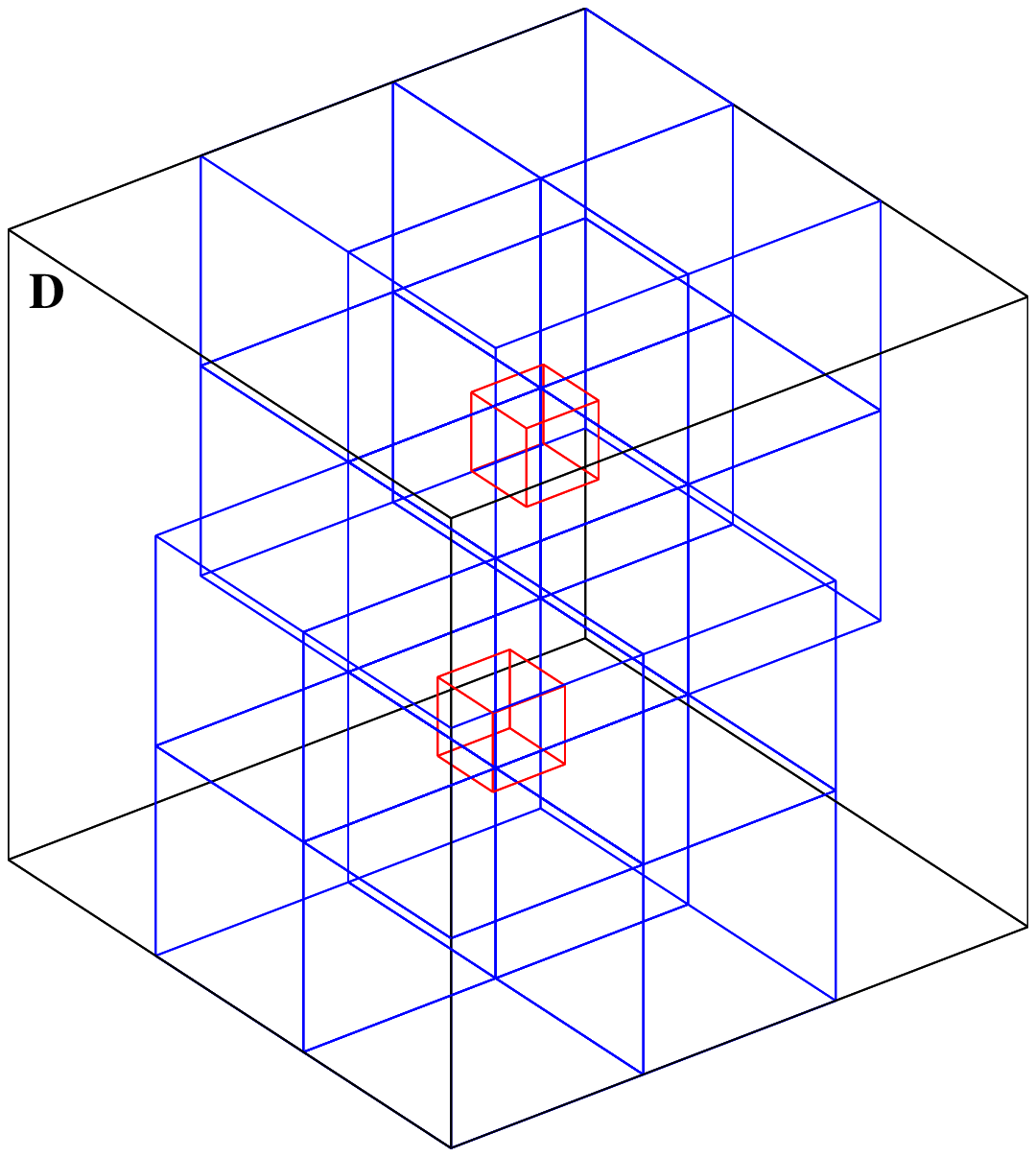}\hfill{}
 \hfill{}\includegraphics[clip,width=0.45\textwidth]{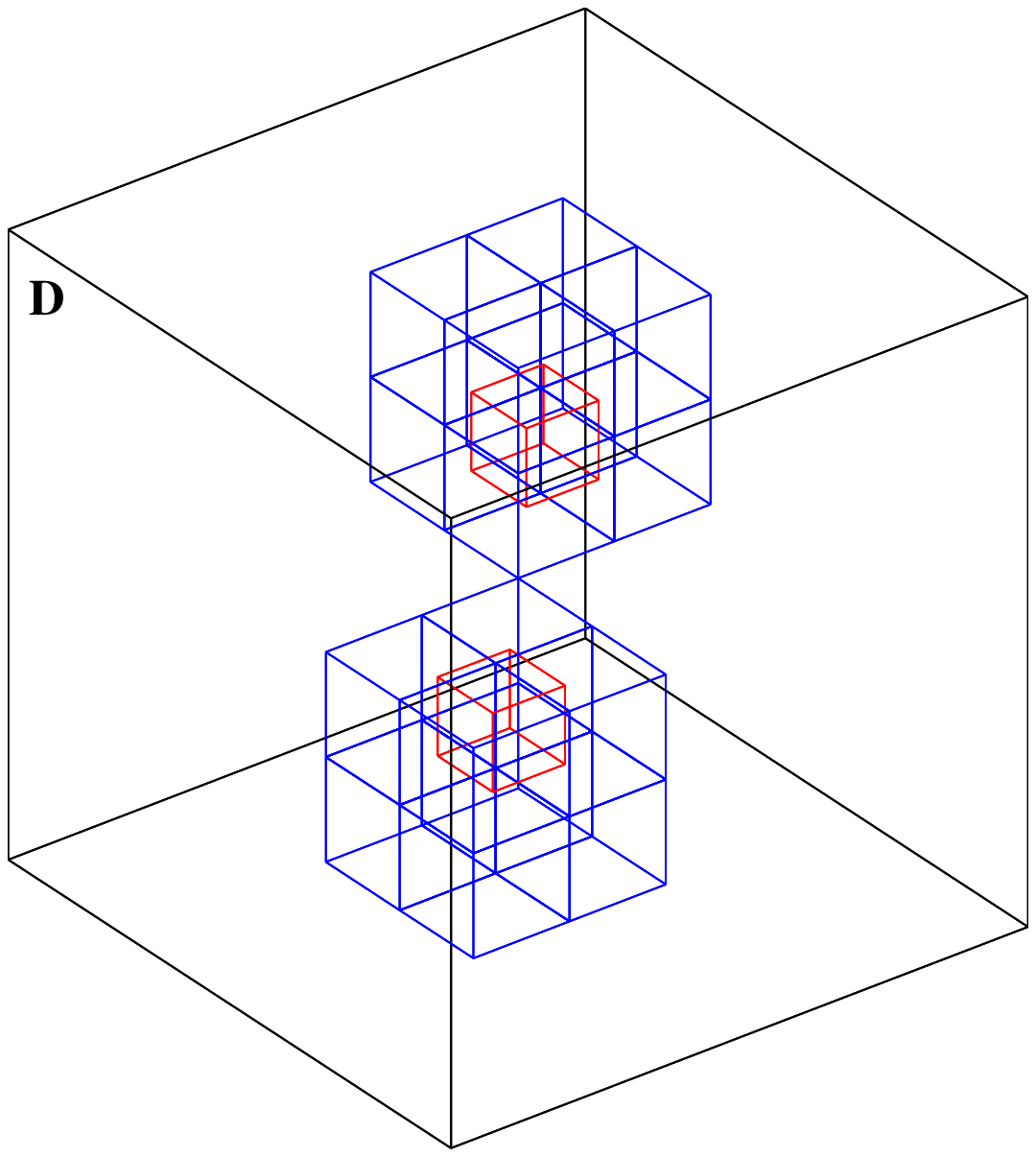}\hfill{}

 \hfill{}(a)\hfill{} \hfill{}(b)\hfill{}

 \hfill{}\includegraphics[clip,width=0.45\textwidth]{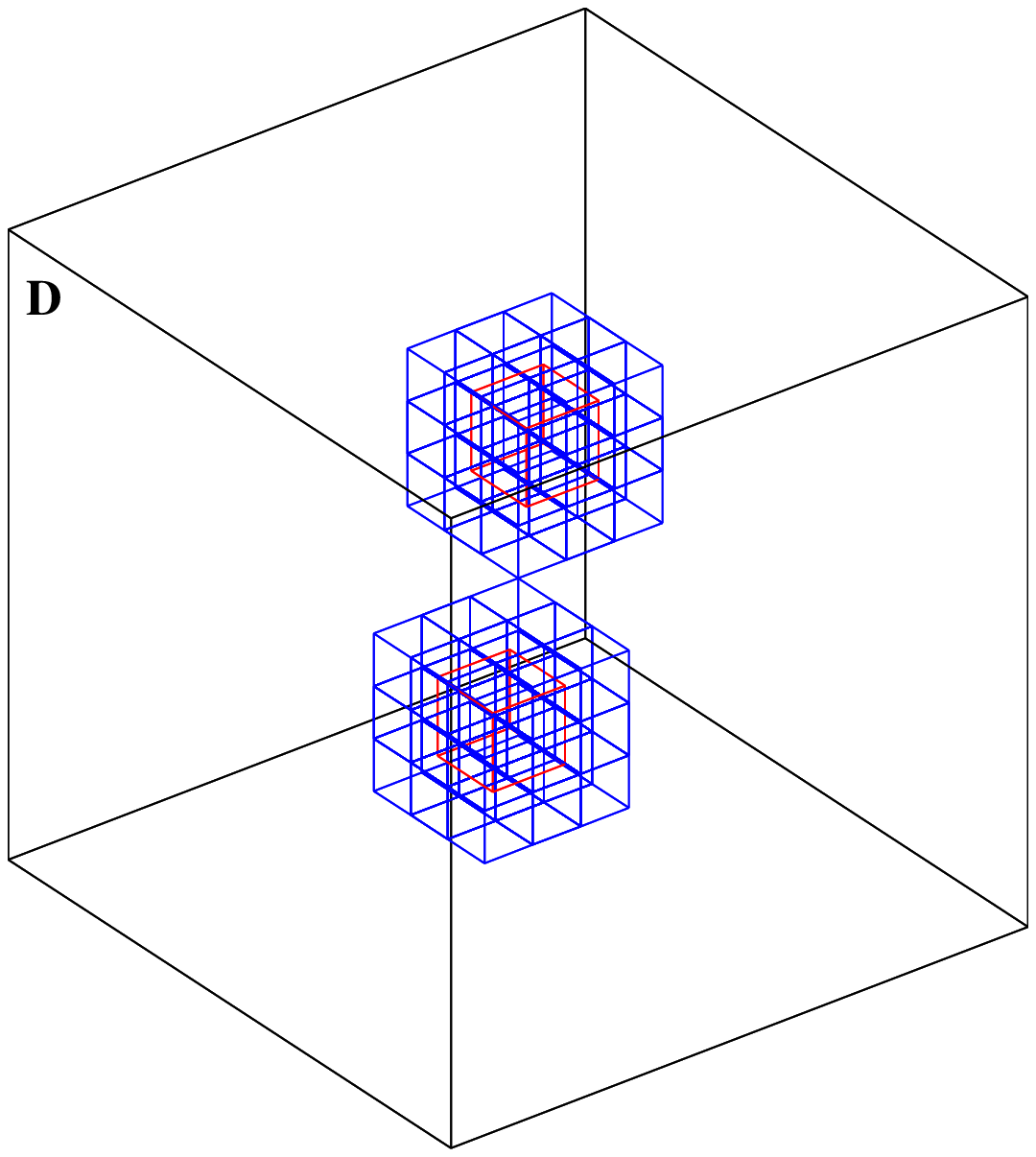}\hfill{}
 \hfill{}\includegraphics[clip,width=0.45\textwidth]{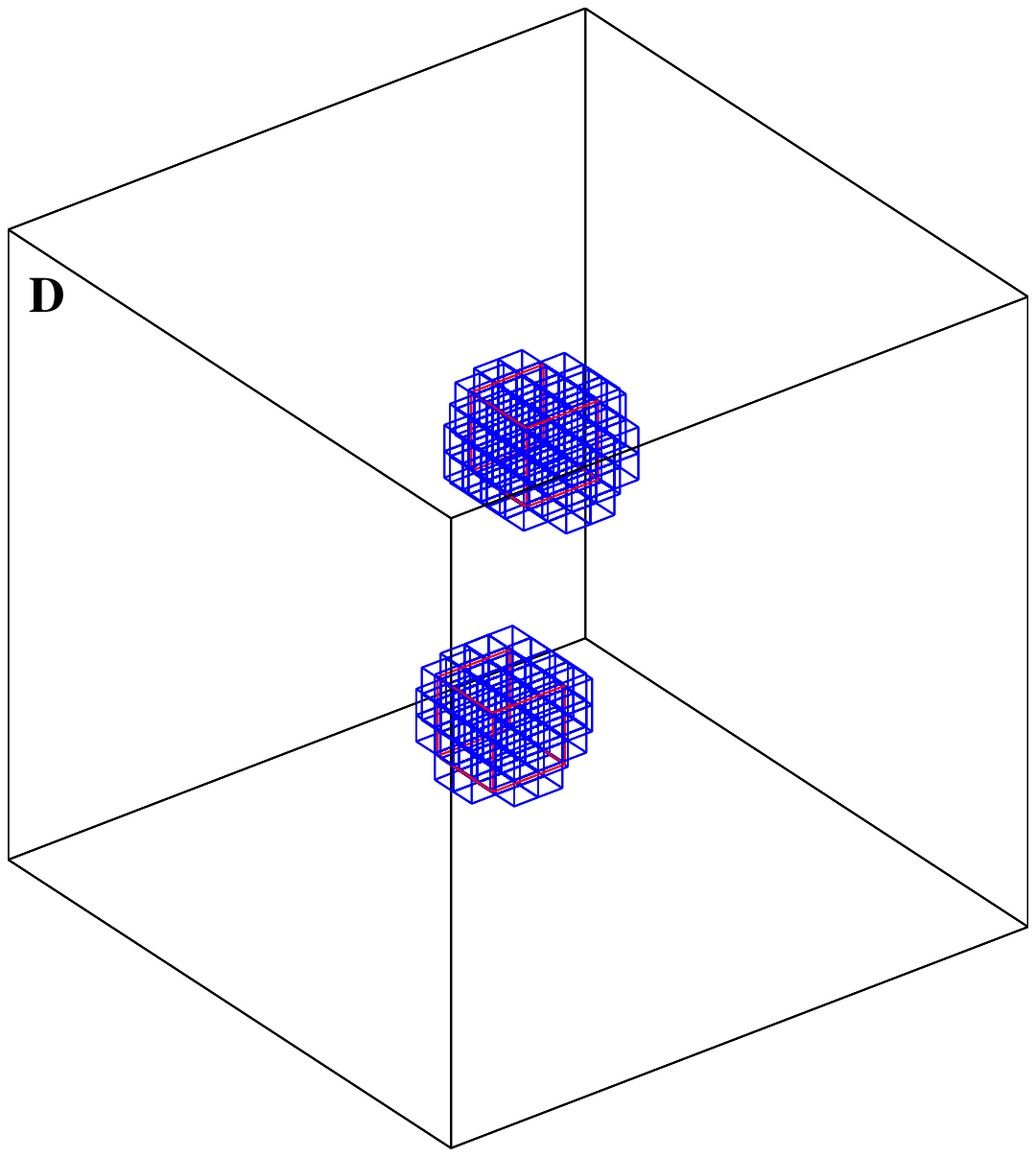}\hfill{}

 \hfill{}(c)\hfill{} \hfill{}(d)\hfill{}

 \caption{\label{fig:simulation4}\emph{Numerical reconstructions by the first 4 iterations}}
 \end{figure}

We take an initial mesh size of $h_0={0.8\lambda}$ in the multilevel algorithm. The mesh refinement during the multilevel algorithm is  carried out 
based on the rule
$h_{k}={0.8\lambda}/{2^{k}}$, where $k$ is the $k$-th refinement, and $h_k$ 
is the mesh size after the $k$th refinement.  
 The numerical reconstructions are shown in Figure \ref{fig:simulation4}.
 Same as for the previous two-dimensional examples, 
 the reconstructions are quite satisfactory and the accurate locations of 
 the scatterer can be achieved, and two inhomogeneous medium objects 
 can be quickly separated.
 
 \begin{figure}[!htb]
 \hfill{}\includegraphics[clip,width=0.45\textwidth]{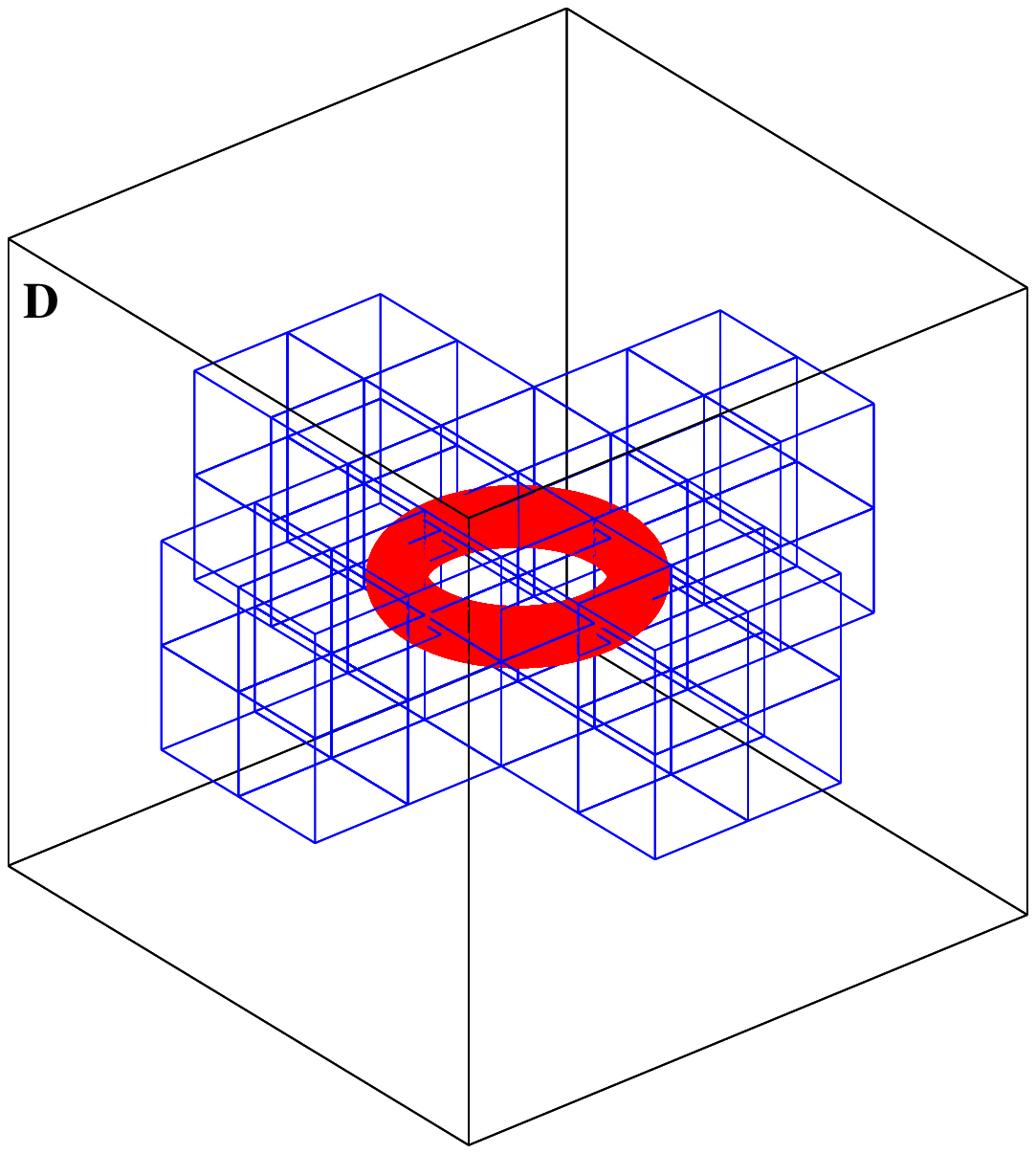}\hfill{}
 \hfill{}\includegraphics[clip,width=0.45\textwidth]{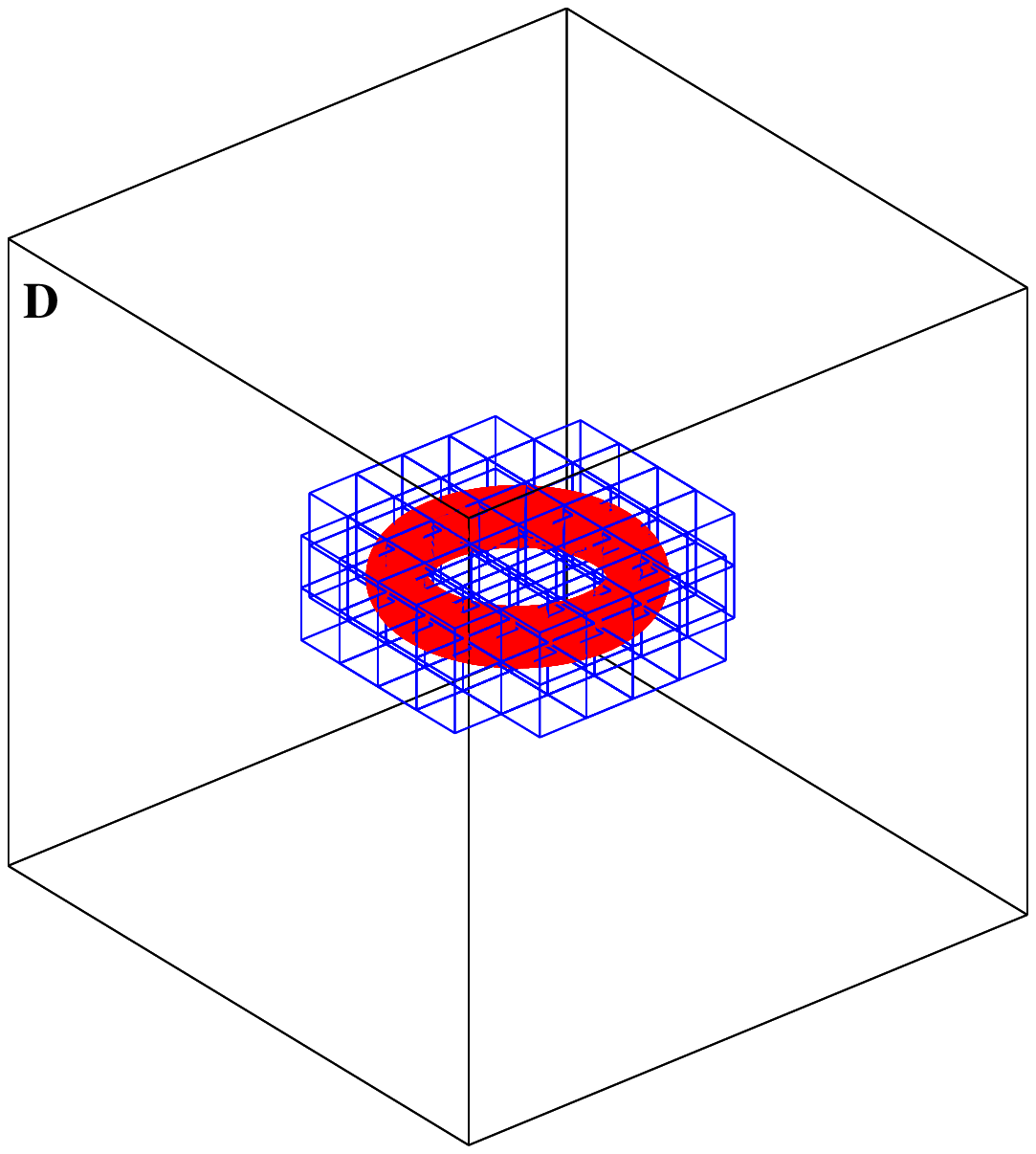}\hfill{}

 \hfill{}(a)\hfill{} \hfill{}(b)\hfill{}

 \hfill{}\includegraphics[clip,width=0.45\textwidth]{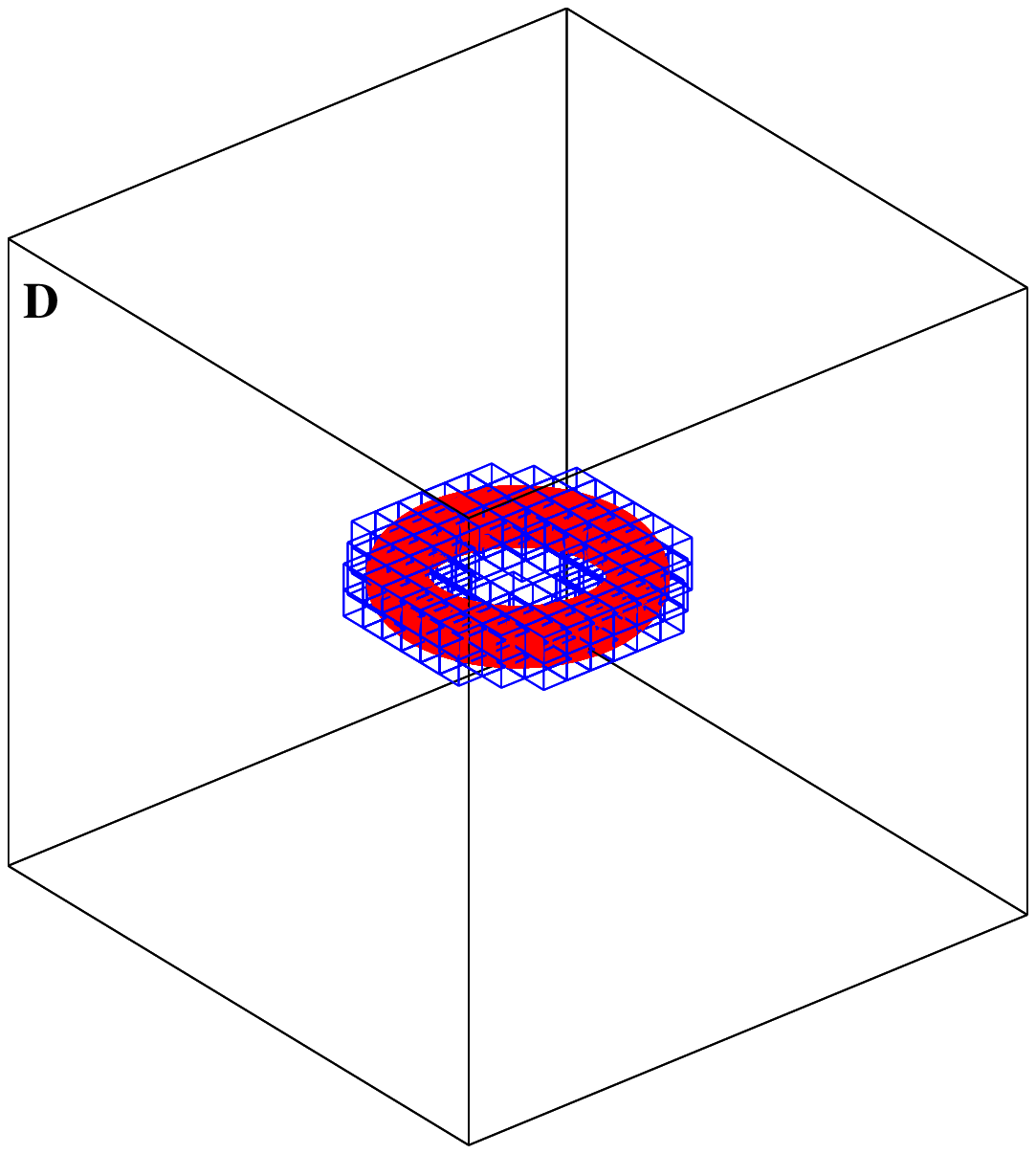}\hfill{}
 \hfill{}\includegraphics[clip,width=0.45\textwidth]{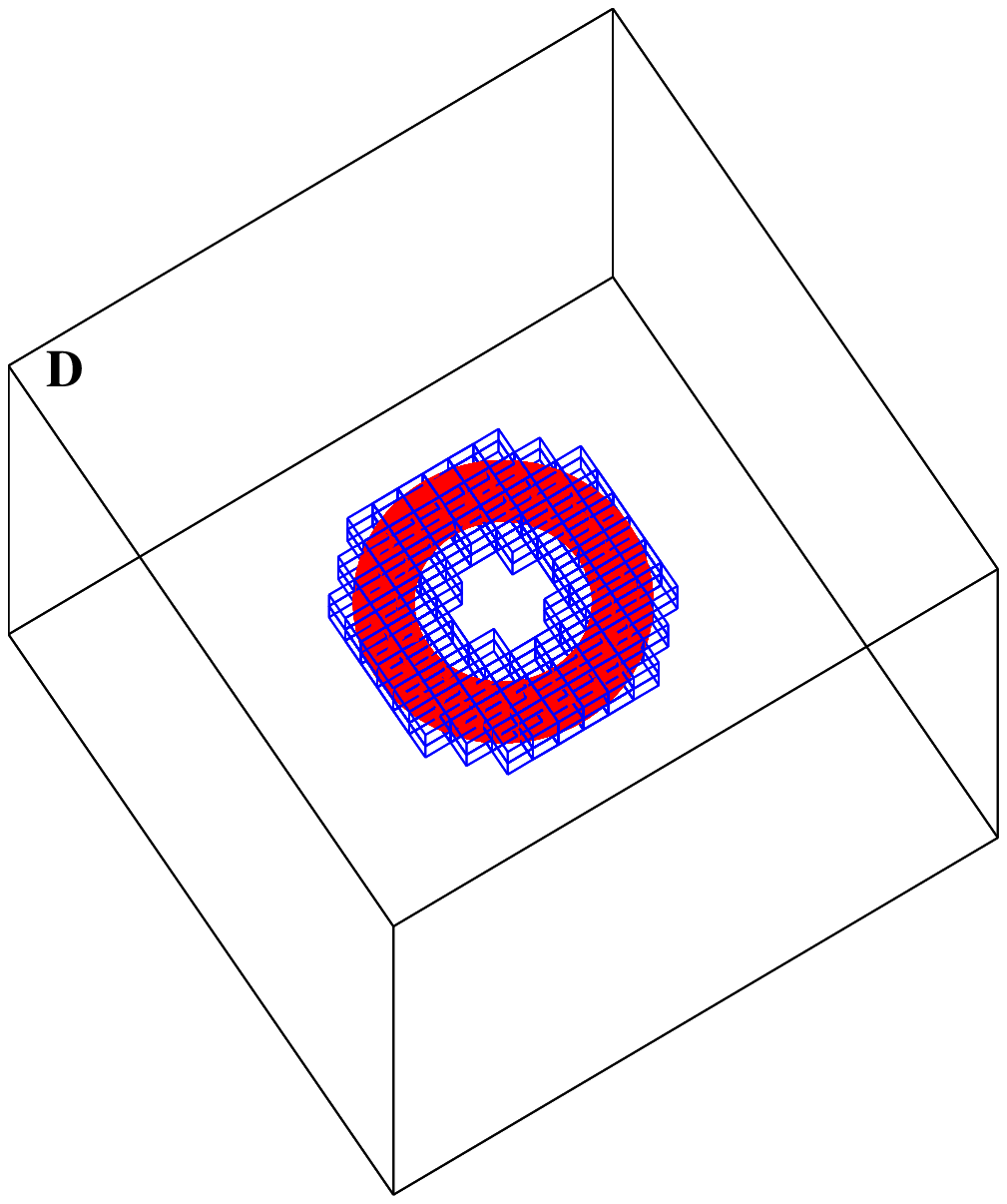}\hfill{}

 \hfill{}(c)\hfill{} \hfill{}(d)\hfill{}

 \caption{\label{fig:simulation5}\emph{Numerical reconstructions by the first 3 iterations; (d) 
 is the same as (c), but viewed in a different angle.}}
 \end{figure}
   
\ss
{\bf Example 5}. In this test we consider a torus scatterer (see Figure~\ref{fig:boxes}(b)), 
with a contrast value $2$. The torus has the following representation,
\[\bigg(R-\sqrt{x^2+y^2}\bigg)^2+z^2=r^2,\] 
where
$r=0.1\lambda$ and $R=0.4\lambda$ ($R$ is the radius from the center of the hole to the center of the torus tube, $r$ is the radius of the tube).
The sampling doman is taken to be 
$\mathbf{D}=[-1.2\lambda,1.2\lambda]^3$, which is about 170 times of the volume 
of the torus. 

We take an initial mesh size of $h_0={0.4\lambda}$ in the multilevel algorithm. The mesh refinement during the multilevel algorithm is  carried out 
based on the rule
$h_{k}={0.4\lambda}/{2^{k}}$, where $k$ is the $k$-th refinement, and $h_k$ 
is the mesh size after the $k$th refinement.  
 The numerical reconstructions are shown in Figure \ref{fig:simulation5}.
 Same as for the previous two-dimensional examples, 
 the reconstructions are quite satisfactory and the accurate locations of 
 the scatterer can be achieved.

\section{Concluding remarks}
 This work proposes a multilevel sampling algorithm 
 which helps locate an initial computational domain  
 for the numerical reconstruction of inhomogeneous media in inverse medium scatterings. 
 The algorithm is an iterative process which starts with a large sampling domain, 
 and reduces the size of the domain iteratively based on 
 the cut-off values, which are computed adaptively by using the updated 
 contrast source strengths and contrast values at each iteration.  
 The iterative algorithm can be viewed actually as  
 a direct method, since it involves only matrix-vector 
operations and does not need any optimization process or 
to solve any large-scale ill-posed linear systems.  
The algorithm works with very few incident fields and its cut-off values are easy 
to compute and insensitive to the sizes and shapes of 
the scatterers, as well as the noise in the data. This is a clear advantage of 
the algorithm over some popular existing sampling methods such as  the linear sampling type methods, 
where the cut-off values are sensitive to the noise and difficult to choose, 
and the number of incident fields can not be small.  
In addition, the multilevel algorithm converges fast and can easily separate 
multiple disjoint scattering components, often with just a few iterations 
    to find a satisfactory initial location of each object. 
    Another nice feature of the new algorithm is that it is self-adaptive, that is, it can
    remedy the possible errors from the previous levels at each current level. 
    With an effective initial location of each object,
    we may then apply any existing efficient but computationally more demanding 
    methods for further refinement of the estimated shape of each 
scattering object as well as for recovery of the contrast profiles of different media.



\begin{thebibliography}{11}
    
     \bibitem{AbuBer} A. Abubakar and P. M. van den Berg, 
     The contrast source inversion method for location and shape reconstructions, 
     Inverse Problems \textbf{18} (2002), pp.\,495-510.
      
    \bibitem{bao05} G. Bao and P. Li,  Inverse medium scattering for the Helmholtz equation at fixed frequency, Inverse Problems \textbf{21} (2005), pp.1621-1641.
   
    \bibitem{BelChaSen} K. Belkebir, P. C. Chaumet and A. Sentenac, Superresolution in total internal reflection tomography, J. Opt. Soc. Amer. A \textbf{22} (2005), pp.\,1889-1897.
    
     \bibitem{Chen09} X. Chen, Application of signal-subspace and optimization methods in reconstructing extended scatterers, J. Opt. Soc. Amer. A \textbf{26} (2009), pp.\,1022-1026.

    \bibitem{Chen10} X. Chen, Subspace-based optimization method for solving inverse-scattering problems, IEEE Trans. Geosci. Remote Sensing \textbf{48} (2010), pp.\,42-49.

    \bibitem{ColKre} D. Colton and R. Kress, Inverse Acoustic and Electromagnetic Scattering Theory, 2nd ed., Springer Verlag, Berlin, 1998. 
    
    \bibitem{Jin12}  K. Ito, B. Jin and J. Zou,  A direct sampling method to an inverse medium scattering 
    problem,  Inverse Problems \textbf{28} (2012), 025003 (11pp). 
    
    
    \bibitem{Kir} A. Kirsch, The MUSIC-algorithm and the factorization method in inverse scattering theory for inhomogeneous media, Inverse Problems \textbf{18} (2002), pp.\,1025-1040.
      
     \bibitem{Lak} A. Lakhtakia, Strong and weak forms of the method of momoents and the coupled dipole method for scattering of time-harmonic electromagnetics fields, Int. J. Modern Phys. C \textbf{3} (1992), pp.\,583-603.
      
    \bibitem{Levy88} B. Levy and C. Esmersoy, Variable background Born inversion by wavefield backpropagation, SIAM J. Appl. Math. \textbf{48} (1988), pp.\,952-972. 
      
    \bibitem{LLZ08} J. Li, H. Liu and J. Zou, Multilevel linear sampling method for inverse scattering problems, SIAM J. Sci. Comp. \textbf{30} (2008), pp.\,1228-1250.
 
    \bibitem{LLZ09} J. Li, H. Liu and J. Zou, Strengthened linear sampling method with a reference ball, SIAM J. Sci. Comp. \textbf{31} (2009), pp.\,4013-4040.
     
     \bibitem{LiuXuZou} K. Liu, Y. Xu and J. Zou, A parallel radial bisection algorithm for inverse scattering problems, Inv. Prob. Sci. Eng. (2012), pp.\,1-13.
        
    \bibitem{MarGruSim} E. A. Marengo, F. K. Gruber and F. Simonetti, Time-reversal MUSIC imaging of extended targets, IEEE Trans. Image Proc. \textbf{16} (2007), pp.\,1967-1984.
    
    \bibitem{Pot} R. Potthast, A survey on sampling and probe methods for inverse problems, Inverse Problems \textbf{22} (2006), pp.\,R1-R47.
    
    \bibitem{BerKle} P. M. van den Berg and R. E. Kleinman, A contrast source inversion method, Inverse Problems \textbf{13} (1997), pp.\,1607-1620.

    \bibitem{BerBro} P. M. van den Berg, A. L. van Broekhoven and A. Abubakar, Extended contrast source inversion, Inverse Problems \textbf{15} (1999), pp.\,1325-1344.

    \end{thebibliography}
\end{document}